\documentclass[11pt,a4paper]{amsart}

\usepackage{amsmath,amstext,amssymb,amsthm,amsfonts}
\usepackage[all,cmtip]{xy}

\theoremstyle{plain}% default
\newtheorem*{theorem*}{Theorem}
\newtheorem*{remark*}{Remark}
\newtheorem*{example*}{Example}

\newtheorem*{conjecture*}{Conjecture}

\theoremstyle{definition}

\theoremstyle{remark}

\newcommand{\bQ}{{\mathbb Q}}
\newcommand{\bR}{{\mathbb R}}

\newcommand{\abs}[1]{\left|{#1}\right|}

\def\quotient#1#2{%
    \raise1ex\hbox{$#1$}\Big/\lower1ex\hbox{$#2$}%
}

\begin{document}

\date{November 18, 2019}
\setcounter{tocdepth}{1}
\title[Tempered distributions and Schwartz functions on definable manifolds]{Tempered distributions and Schwartz functions on definable manifolds}
 \author{Ary Shaviv}

\email{ ashaviv@math.princeton.edu}
%\thanks{$^\dag$Supported in part by ERC StG grant 637912 and ISF grant 249/17}

\maketitle

\begin{abstract}
We define the spaces of Schwartz functions, tempered functions and tempered
distributions on manifolds definable in polynomially bounded o-minimal structures. We show that all the classical properties
that these spaces have in the Nash category, as first studied in Fokko du
Cloux's work, also hold in this generalized setting. We also show that on manifolds definable in o-minimal structures that are not polynomially bounded, such a theory can not be constructed. We present some
possible applications, mainly in representation theory.
\end{abstract}

\tableofcontents

\section{Introduction}\label{sec-intro}

The theory of distributions, that was first introduced
by Laurant Schwartz in \cite{schwartz}, has great importance in the rigourous
definition of Fourier transform, in a manner that gives a meaningful insight
on Heisenberg's uncertainty principle. A natural class of distributions to
consider is the class of tempered distributions -- these are functionals
on the space of compactly supported smooth functions, that extend to the
space of Schwartz functions. A smooth real valued function on $\bR^n$ is
called Schwartz if it remains bounded on all of $\bR^n$, even after an arbitrary
algebraic differential operator is applied to it. In \cite{dC}, Fokko du
Cloux generalized the notion of Schwartz functions to an arbitrary affine
Nash manifold -- loosely speaking, an affine Nash manifold is a closed semi-algebraic
subset of $\bR^n$ who is also a $C^\infty$-smooth sub-manifold. His motivation
came from representation theory: let $G$ be a Nash group (a semi-algebraic
Lie group), let $H$ be a closed subgroup of $G$, and let $\pi$ be an arbitrary
representation of $H$. The notion of Schwartz functions on affine Nash manifolds
enables in some cases the definition of Schwartz induction of representations, namely the
construction of a new representation $\mathcal{S}\text{Ind}_H^G(\pi)$ of
$G$. In \cite{AG}, Avraham Aizenbud and Dmitry Gourevitch generalized the
theory of Schwartz functions to the category of general (not necessarily
affine) Nash manifolds, and in \cite{ES}, Boaz Elazar and the author constructed
a theory of Schwartz functions on (not necessarily smooth) affine algebraic
varieties. An important object in the above theories is the space of tempered
functions -- a smooth function on an open subset of $\bR^n$ is tempered if
it, together with all of its derivatives, goes to infinity as approaching
the boundary (or infinity), not faster than one over some polynomial in the
distance to the boundary. A wide study of these spaces of functions was done
by Masaki Kashiwara and Pierre Schapira in \cite{KS1,KS2}. Although Bernard
Malgrange showed in \cite{Ma} that these spaces of functions do not form
a sheaf in general, Kashiwara and Schapira essentially showed that "sheaf
property" does hold whenever $\L$ojasiewicz's inequality holds.

\

Nash manifolds are the simplest example of manifolds definable in a o-minimal structures (the reader not accustomed to the  model-theoretic notion of "definability in structures" will find an introductory review in subsection \ref{intro-to-definablility}). In the model-theoretic language, Nash manifolds are simply manifolds definable in the ordered field of real numbers $\bR$, i.e. in the structure $(\bR,<,0,1,+,\cdot)$. This structure is polynomially bounded -- any function on $\bR$ that is definable in it is bounded by some polynomial. Polynomially bounded structures have "tame geometry", in the sense that $\L$ojasiewicz's inequality, with its many implications, holds.

\

In this paper we construct the theory of Schwartz functions, tempered functions
and tempered distributions on manifolds that are definable in arbitrary polynomially bounded o-minimal structures. We show that their "tame geometry" implies that all
the classical properties these spaces have in the Nash case also hold in
the this framework. For manifolds definable in o-minimal structures that are not polynomially bounded, we show that the theory of Schwartz functions is always ill-defined (this follows from Miller's Dichotomy Theorem saying that the exponential function is always definable in this case) -- in this sense our result in ultimate.

\

{\bf Motivation.} The first motivation for our theory is realizing du Cloux's
Schwartz induction as a space of sections: using the notation above, one
often thinks of induced representation as a space of sections of the vector
bundle $(G\times \pi)/_\sim\to G/H$. Realizing Schwartz induction in this
language requires the definition of Schwartz sections. However, even when
$G\times \pi$ is Nash, the quotient $(G\times \pi)/_{\sim}$ is usually not
Nash, and so Schwartz sections are not defined neither in du Cloux's setting,
nor in Aizenbud-Gourevitch's. The results in this paper enable the notion of Schwartz sections in some of these cases. This is discussed in
more detail in Appendix \ref{appendix-rep-theory}, together with another
motivation from representation theory.

\

A second motivation is the following theorem characterizing Schwartz spaces
on open subsets in the Nash category (see \cite[Theorem 3.23]{ES} for a similar
result in the affine algebraic case):

\subsection{Theorem (\cite[Theorem 5.4.1]{AG})}\label{AG_motivating_theorem}Let
$M$ be a Nash manifold and $U\subset M$ be an open (semi-algebraic) subset.
Let $\mathcal{S}(M)$ and $\mathcal{S}(U)$ be the spaces of Schwartz functions
on $M$ and on $U$ respectively, and $W\subset\mathcal{S}(M)$ be the closed
subspace defined by: $$W:=\{\phi\in\mathcal{S}(M)|\phi\text{ vanishes with
all its derivatives on }M\setminus U\}.$$ Then, restriction to $U$ and extension
by zero give an isomorphism $W\cong\mathcal{S}(U)$.

\

{\bf The basic ideas.} Inspired by Theorem \ref{AG_motivating_theorem}, it
is natural to define Schwartz functions on an arbitrary open, not necessarily
semi-algebraic, $U\subset \bR^n$: these should be all real valued functions
on $U$, whose extensions by zero to all of $\bR^n$ are Schwartz functions
on $\bR^n$ that vanish with all their derivatives on $\bR^n\setminus U$.
This definition (Definition \ref{schwartz space for euclidean open definition}
below) is indeed the starting point of this paper.

\

Affine definable manifolds in a polynomially bounded o-minimal structure are manifolds that may be realized "in a definable way" as open subsets of $\bR^n$. On such manifolds Schwartz functions are naturally defined
by considering any open definable embedding in $\bR^n$. On general definable manifolds in a polynomially bounded o-minimal structure we define Schwartz functions to be finite sums of extensions by
zero of Schwartz functions on open affine definable sub-manifolds. We
denote the Fr\'echet space of Schwartz functions on such a definable manifold
$M$ by $\mathcal{S}(M)$. We also study the following two spaces: the space
of tempered distributions on $M$  -- denoted by $\mathcal{S}^*(M)$
-- the space of continuous linear functionals on the Schwartz space; the
space of tempered functions -- denoted by $\mathcal{T}(M)$ -- the space of
all functions that have the Schwartz space as a module under point-wise multiplication,
i.e. a real valued function $t$ is tempered if for any $s\in\mathcal{S}(M)$,
also $t\cdot s\in\mathcal{S}(M)$. In particular we show that in the affine
relatively compact case this notion of tempered functions and Kashiwara and Schapira's notion
coincide.

\

Recall that a pre-sheaf of functions $\mathcal{G}$ on a topological space
$M$ is called a sheaf if for any open subset $U\subset M$ and any open cover
$U=\bigcup\limits_{\alpha\in A}U_\alpha$, the following sequence is exact:

\begin{center}
\leavevmode
\xymatrix{
0 \ar[rr] & & \mathcal{G}(U) \ar[rr]^{\text{Res}_1} & & \prod\limits_{\alpha\in
A}\mathcal{G}(U_\alpha) \ar[rr]^(.4){\text{Res}_2} & & \prod\limits_{\alpha\in
A}\prod\limits_{\beta\in A}\mathcal{G}(U_\alpha\cap U_\beta),}
\end{center}
where the $\alpha$-coordinate of $\text{Res}_1(f)$ is $f|_{U_\alpha}$, and
the $(\alpha,\beta)$-coordinate of $\text{Res}_2(\prod\limits_{\alpha\in
A}f_\alpha)$ is $f_\alpha|_{U_\alpha\cap U_\beta}-f_\beta|_{U_\alpha\cap
U_\beta}$. Definable manifolds naturally form restricted topological spaces (and not topological spaces), and so it is natural to define a sheaf on them as a pre-sheaf such that the sequence above is exact
in all cases where the index set $A$ is finite, and $U_\alpha$ is definable
for any $\alpha\in A$ (this can be formulated by defining a Grothendieck
topology on $M$, but we intentionally avoid this language). Similarly we
define other sheaves and cosheaves on definable manifolds.

\

{\bf Main results.} We show that all the classical properties the spaces
of Schwartz functions, tempered functions and tempered distributions have
in the Nash and the algebraic cases also hold in our settings. In particular
we prove the following:

\

Let $M$ be a smooth manifold definable in some polynomially bounded o-minimal structure and $M=\bigcup\limits_{i=1}^{m}U_i$ be
a finite definable open cover, then:
\begin{enumerate}
  \item {\bf Tempered partition of unity on definable
manifolds} (Proposition
\ref{tempered-partition-of-unity-definable} and Lemma \ref{lemma_temp_move_s_to_s_definable}).
There exist non-negative tempered functions $\varphi_1,\dots,\varphi_m\in\mathcal{T}(M)$
such that the following hold:\begin{enumerate}
                                                       \item $\text{supp}(\varphi_i)\subset
U_i$ for any $i=1,\dots,m$;
                                                       \item $\sum\limits_{i=1}^{m}\varphi_i(x)=1$
for any $x\in M$;
                                                       \item $(\varphi_i\cdot
s)|_{U_i}\in\mathcal{S}(U_i)$ for any $s\in\mathcal{S}(M)$ and for any $i=1,\dots,m$.
                                                     \end{enumerate}
                                                       \item {\bf Sheaf property of tempered functions on definable
manifolds}
(Proposition \ref{tempered-is-a-sheaf-definable}). The assignment of the
space of tempered functions to any definable open $U\subset M$, together
with the restriction of functions, form a sheaf.

  \item {\bf Sheaf property of tempered distributions on definable
manifolds}
(Proposition \ref{tempered-distributions-is-a-sheaf-definable}). The assignment
of the space of tempered distributions to any definable open $U\subset
M$, together with the restriction of functionals, form a flabby sheaf.
  \item {\bf Cosheaf property of Schwartz functions on definable manifolds}
(Proposition \ref{schwartz_functions-is-a-cosheaf-definable}). The assignment
of the space of Schwartz functions to any definable open $U\subset M$,
together with extension by zero, form a flabby cosheaf.
  \item {\bf Compactly supported functions} (Theorem \ref{prop_any_compact_supp_fun_is_schwartz}).
$C^\infty_c(M)\subset\mathcal{S}(M)$. In particular if $M$ is compact, then
$\mathcal{S}(M)=C^\infty(M)$.
\item {\bf Schwartz functions with full support} (Corollary \ref{thm-schwartz-with-full-support_cor}).
There exists $s\in\mathcal{S}(M)$
such that $s(x)>0$ for any $x\in M$.
\end{enumerate}

{\bf Structure of this paper.} In {\bf Section \ref{sec-prelim}} all
preliminary results used in this paper are presented. These include mainly a quick review on definable manifolds in o-minimal structures, and some
classical results on polynomially bounded o-minimal structures, exponential structures, Fr\'echet
spaces, and finally on the classical theory of Schwartz functions on $\bR^n$.

In {\bf Section \ref{sec-arbitrary}} we define and study the
spaces of Schwartz functions, tempered functions and of tempered distributions
on arbitrary open subsets of $\bR^n$.

In {\bf Section \ref{sec-embedded}} we define and study
the spaces above in the special case of open subsets that are definable in some polynomially bounded o-minimal structure, and prove properties
(1)-(4) in this case.

In {\bf Section \ref{sec-affine}} we
define and study the spaces of Schwartz functions, tempered functions and
of tempered distributions on affine manifolds that are definable in some polynomially bounded o-minimal structure. We show that
properties (1)-(4) hold in this case, as an easy implication of Section \ref{sec-embedded}.

In {\bf Section \ref{sec-general}} we
extend the results of Section \ref{sec-affine}
to general (non-affine) manifolds that are definable in some
polynomially bounded o-minimal structure, and also prove property (5).

{\bf Section \ref{sec-whitney}} is devoted to proving property (6): on
any definable manifold in a polynomially bounded o-minimal structure there exists a strictly positive Schwartz function.

In {\bf Section \ref{sec-miller}} we use Miller's Dichotomy Theorem in order to show that on manifolds definable in o-minimal structures that
are not polynomially bounded, the theory of Schwartz functions
is always ill-defined. 

In {\bf Appendix \ref{appendix-rep-theory}} we explain how to apply
our theory in representation theory.

In {\bf Appendix \ref{appendix-conjectures}} we present a conjecture on $C^\infty$-diffeomorphisms vs. definable $C^\infty$-diffeomorphisms.

In {\bf Appendix \ref{appendix-fractals}} we suggest a possible invariant
in our theory that might be used in fractal geometry.

Finally, in {\bf Appendix \ref{appendix-nash}} we discuss
the relations between our theory and the theory of Schwartz functions on
Nash manifolds. In particular we show that our theory generalizes the latter and deduce a useful corollary.

\

{\bf Conventions and Notations.} Throughout this paper the base field
is $\mathbb{R}$ -- the real numbers. By the notion \emph{smooth} we always
mean infinitely smooth, e.g. a smooth function on the real line lies in $C^\infty(\mathbb{R})$,
i.e. it is continuously differentiable $k$ times, for any $k\in\mathbb{N}$.
When we write \emph{diffeomorphism} we always mean smooth diffeomorphism.
\emph{Manifolds} are always smooth manifolds. We set $\mathbb{N}_0:=\mathbb{N}\cup\{0\}$,
i.e. $\mathbb{N}_0$ is the set of all non-negative integers. If $U\subset\bR^n$
is an open subset, $f:U\to\bR$ is a smooth function and $\alpha_1,\dots,\alpha_n\in\mathbb{N}_0$,
we use multi-index notation for derivatives: $\alpha:=(\alpha_1,\alpha_2,\dots,\alpha_n)$,
$\abs{\alpha}:=\sum\limits_{i=1}^{n}\alpha_i$, $f^{(\alpha)}=\frac{\partial^{|\alpha|}f}{\partial
x_1^{\alpha_1}\partial x_2^{\alpha_2}\cdot\cdot\cdot\partial x_n^{\alpha_n}}$
if $\abs{\alpha}\neq0$, and $f^{(\alpha)}=f$ if $\abs{\alpha}=0$. For any
$x\in U$ we denote $x^\alpha:=\prod\limits_{i=1}^{n}x_i^{\alpha_i}$, where
$x_i$ is the $i^{th}$ coordinate of $x$, and denote by $|x|$ the standard
Euclidean norm of $x$. For two subsets $Z_1,Z_2\subset U$ we set $\text{dist}(Z_1,Z_2):=\inf\limits_{z_1\in
Z_1,z_2\in Z_2}\abs{z_1-z_2}$. We also set $\text{dist}(x,Z_2):=\text{dist}(\{x\},Z_2)$.
When $X$ is any set, and $Y\subset X$ is any subset, we denote by $\text{Ext}_Y^X$
the "extension by zero" operator that takes a real valued function on $Y$
and returns a real valued function on $X$: for any $f:Y\to\mathbb{R}$, $\text{Ext}_Y^X(f)(x):=f(x)$
for all $x\in Y$, and $\text{Ext}_Y^X(f)(x):=0$ for all $x\in X\setminus
Y$. We denote by $\text{Res}_Y^X$ the natural restriction operator that takes
a real valued function on $X$ and returns a real valued function on $Y$:
for any $g:X\to\mathbb{R}$, $\text{Res}_Y^X(g)(y):=g(y)$ for all $y\in Y$.

\

{\bf Acknowledgments.} The author is grateful to Avraham Aizenbud, Saugata Basu, Itai Benjamini ,Edward Bierstone, Gal Binyamini, Charles Fefferman, Andrei Gabrielov, Dmitry Gourevitch, Bo'az Klartag, Pierre D. Milman, Dmitry Novikov, Ya'akov Peterzil, Pierre Schapira and Patrick Speissegger for many valuable discussions. A special thanks goes to Yotam Hendel for providing some helpful comments on a very early version of this paper. Finally, the author wishes to express his gratitude to the anonymous referee, for carefully reading the manuscript and providing very insightful and helpful comments.

In the final stages of this research the author was supported by AFOSR Grant FA9550-18-1-069.  The early stages of this research were done while the author was a Ph.D. student
in Weizmann Institute of Science, then the author was supported in part by
ERC StG grant 637912 and ISF grant 249/17.

\section{Preliminaries}\label{sec-prelim}

\subsection{Definable manifolds in o-minimal structures}\label{intro-to-definablility}Let $\mathcal{R}:=(\bR,<,0,1,+,\cdot,\dots)$ be an expansion of the ordered field of real numbers. We also call $\mathcal{R}$ a \emph{structure}. A subset of $\bR^n$ ($n\in\mathbb{N}$) is said to be \emph{definable in $\mathcal{R}$} (or simply \emph{definable} if the structure is clear from the context) if it is first-order definable in $\mathcal{R}$ with parameters from $\bR$. A
map $f:X\to Y$ ($X\subset\bR^n,Y\subset\bR^m,n,m\in\mathbb{N}$) is said to
be definable if its graph is definable in $\bR^{n+m}$. A structure is called \emph{o-minimal} (o stands for \emph{order}) if the definable subsets of $\bR$ are exactly all finite unions of intervals (including intervals of infinite length and singletons). Two common examples of such structures are the semi-algebraic category that consists of sets and maps definable in $(\bR,<,0,1,+,\cdot)$, and the globally subanalytic category that consists of sets and maps definable in $\bR_{\text{an}}$ (see \cite[2.5]{vdDM} for the exact definition of $\bR_{\text{an}}$ and for more examples; also see \ref{example-good-structure} below for another detailed example -- $\bR_{\text{an}}^\bR$).

We use similar terminology to that of \cite{PS}: a \emph{definable smooth atlas on a set $X$ relative to an o-minimal structure $\mathcal{R}$} is a finite family of pairs $\{(U_i,g_i)\}_{i\in I}$ (called \emph{definable charts}, or simply \emph{charts}) such that for any $i,j\in I$: $U_i$ is a subset of $X$, $g_i:U_i\to\bR^{n(i)}$ is an injection, $g_i(U_i)\subset\bR^{n(i)}$ is a definable smooth sub-manifold, $g_i(U_i\cap U_j)$ is definable and open in $g_i(U_i)$, $g_j\circ g_i^{-1}:g_i(U_i\cap U_j)\to g_j(U_i\cap U_j)$ is a definable smooth diffeomorphism, and finally $X=\bigcup\limits_{i\in I}U_i$. Two definable smooth atlases on a set $X$ relative to an o-minimal structure $\mathcal{R}$ are called \emph{equivalent} is their union is again a definable smooth atlas. A \emph{definable smooth manifold of dimension $n$} (or simply a definable manifold) is a set $M$ together with an equivalence class of definable smooth atlases, such that one of the atlases in this equivalence class has the form $\{(U_i,g_i)\}_{i\in I}$ with each $g_i(U_i)$ is an open subset of $\bR^n$. Note that on such a definable manifold there exists a unique Euclidean topology induced by these atlases, and so it also defines a unique topological space. Moreover, it uniquely defines a smooth manifold in the usual (not "definable") sense.  Also note that any definable open subset $U\subset\bR^n$ naturally defines a definable manifold (by taking the equivalence class of the trivial atlas $\{(U,Id)\}$). 

A subset $M'\subset M$ of a definable manifold
is called a \emph{definable subset} if there exists an atlas of $M$ of
the form $\{(U_i,g_i)\}_{i\in
I}$ such that the sets $g_i(M'\cap U_i)$ are all definable (and thus all
atlases of $M$ have this property). Two definable manifolds $M_1,M_2$ of dimensions $m_1,m_2$ uniquely define a definable manifold $M_1\times M_2$ by the following procedure: say $M_1$ has an atlas $\{(U_i,g_i)\}_{i\in
I}$ and $M_2$
has an atlas $\{(V_j,h_j)\}_{j\in
J}$, then a defining atlas for $M_1\times M_2$ is $\{(U_i\times V_{j},(g_i,h_j))\}_{i\in
I,j\in J}$. This construction is independent of the choices of atlases made. A map from $M_1$ to $M_2$ is called a \emph{definable map} if its graph is a definable subset in $M_1\times M_2$. In particular, a \emph{definable smooth map} is a map that is both smooth (when we think of $M_1$ and $M_2$ as usual smooth manifolds) and definable (when we think of $M_1$ and $M_2$ as definable smooth
manifolds).
\subsubsection{Definable restricted topology} A \emph{restricted topology} on a set $X$ is a collection $\sigma$ of subsets of $X$ that is closed under finite intersections and \emph{finite} unions (this is sometimes called a Grothendieck topology, but we intentionally avoid this terminology). Subsets in $\sigma$ are called \emph{open} (with respect to the restricted topology $\sigma$), and subsets whose complements are in $\sigma$ are called closed. When we say a restricted topological space $X$ we always bear in mind that $\sigma$ is implicitly also given. 

Fix a definable manifold $M$. The collection of all definable open subsets of $M$ does not form a topology, as infinite unions of such sets is usually not a definable subset of $M$. This collection, however, defines the unique \emph{definable restricted topology on $M$}. Note that in general the notion of closure is not well defined in restricted
topological spaces, but it is always well defined in the case of definable restricted topology of a definable manifold (this follows immediately from the fact that the notion of "closure" in $\bR^n$ is first-order definable in $\mathcal{R}$
with parameters from $\bR$). The following Lemma will be very useful to us (as it is not completely standard a proof is given): 

\subsubsection{Lemma (normality of definable restricted topology)}\label{normality_of_definable_lemma}Let $M$ be a definable manifold in an o-minimal structure $\mathcal{R}$, and let $X,Y\subset M$ be definable disjoint closed subsets. Then, there exist definable disjoint open
subsets $\tilde{X},\tilde{Y}\subset M$, such that $X\subset\tilde{X}$
and $Y\subset\tilde{Y}$.

\

Proof: It is clearly enough to prove for the case where $M$ has an atlas consisting of one chart (later we will call such a manifold "affine"), and then we may assume $M\subset\bR^n$ is an open definable subset. Now set $$\tilde{X}=\{u\in M:\exists x\in X:\abs{x-u}<\frac{\text{dist}(x,Y)}{10}\};$$
$$\tilde{Y}=\{u\in M:\exists y\in Y:\abs{y-u}<\frac{\text{dist}(y,X)}{10}\}.$$\qed

\subsubsection{Sheaves and cosheaves on restricted topological spaces} A \emph{pre-sheaf} $\mathcal{F}$ on a restricted topological space $X$ is an assignment $U\mapsto\mathcal{F}(U)$ for any $U\in\sigma$ of an Abelian group, vector space etc., and for any inclusion $V\subset U$ with $V,U\in\sigma$ a "restriction" morphism (in the appropriate category) $res^{U}_V:\mathcal{F}(U)\to\mathcal{F}(V)$ such that $res^U_U=Id$ for all $U\in \sigma$, and for any $W\subset V\subset U$ (all in $\sigma$), $res^V_W\circ res^U_V=res^U_W$. A \emph{sheaf} $\mathcal{F}$ on a restricted topological space $X$ is a pre-sheaf such that the following sequence is exact for any collection $U_1,\dots U_n\in\sigma$: 

\begin{center}
\leavevmode
\xymatrix{
0 \ar[rr] & & \mathcal{F}(U) \ar[rr]^{\text{Res}_1} & & \prod\limits_{i=1}^m\mathcal{F}(U_i) \ar[rr]^(.4){\text{Res}_2} & & \prod\limits_{i=1}^{m-1}\prod\limits_{j=i+1}^m\mathcal{F}(U_i\cap U_j),}
\end{center}
where $U:=\bigcup\limits_{i=1}^m U_i$, the map $\text{Res}_1$ is defined by $\text{Res}_1(\xi):=\prod\limits_{i=1}^m res^U_{U_i}(\xi)$, and the map $\text{Res}_2$ is defined
by $\text{Res}_2(\prod\limits_{i=1}^m \xi_i):=\prod\limits_{i=1}^{m-1}\prod\limits_{j=i+1}^m res^{U_i}_{U_i\cap U_j}(\xi_i)-res^{U_j}_{U_i\cap U_j}(\xi_j)$.
A sheaf is called \emph{flabby} is all the restriction morphisms are onto.

A \emph{pre-cosheaf} $\mathcal{F}$ on a restricted topological space $X$ is
an assignment $U\mapsto\mathcal{F}(U)$ for any $U\in\sigma$ of an Abelian
group, vector space etc., and for any inclusion $V\subset U$ with $V,U\in\sigma$
an "extension" morphism (in the appropriate category) $ext^{U}_V:\mathcal{F}(V)\to\mathcal{F}(U)$
such that $ext^U_U=Id$ for all $U\in \sigma$, and for any $W\subset V\subset
U$ (all in $\sigma$), $ext^U_V\circ ext^V_W=ext^U_W$. A \emph{cosheaf} $\mathcal{F}$
on a restricted topological space $X$ is a pre-cosheaf such that the following
sequence is exact for any collection $U_1,\dots U_n\in\sigma$:

\begin{center}
\leavevmode
\xymatrix{
\bigoplus\limits_{i=1}^{m-1}\bigoplus\limits_{j=i+1}^{m}\mathcal{F}(U_i\cap
U_j) \ar[rr]^(.6){\text{Ext}_1} & & \bigoplus\limits_{i=1}^{m}\mathcal{F}(U_i)
\ar[rr]^{\text{Ext}_2} & & \mathcal{F}(U) \ar[rr] & & 0,}
\end{center}
where $U:=\bigcup\limits_{i=1}^m U_i$, the $k^{th}$ coordinate of $\text{Ext}_1(\bigoplus\limits_{i=1}^{m-1}\bigoplus\limits_{j=i+1}^{m}\xi_{ij})$
is defined by $\sum\limits_{i=1}^{k-1} \text{ext}_{U_k\cap U_i}^{U_k}(\xi_{ik})-\sum\limits_{i=k+1}^{m}
\text{ext}_{U_k\cap U_i}^{U_k}(\xi_{ki})$, and the map $Ext_2$ is defined by $\text{Ext}_2(\bigoplus\limits_{i=1}^{m}\xi_i):=\sum\limits_{i=1}^{m}\text{ext}^{U}_{U_i}(\xi_i)$. A cosheaf is called \emph{flabby} is all the extension morphisms are injective.

\subsection{Polynomially bounded structures}A structure $\mathcal{R}$ is called \emph{polynomially bounded} if for any definable function $f:\bR\to\bR$ there exist $N\in \mathbb{N}$ and $M\in\bR_{>0}$ such that $\abs{f(x)}\leq x^N$ whenever $x>M$.

\subsubsection{Theorem ($\L$ojasiewicz's inequality, \cite[4.14]{vdDM})}\label{Lojasiewicz_inequality} Let $\mathcal{R}$ be a polynomially bounded o-minimal structure, let $K\subset\bR^n$ be a compact subset definable in $\mathcal{R}$, and let $f,g:K\to\bR$ be continuous functions definable in $\mathcal{R}$ such that $\{x\in K:f(x)=0\}\subset\{x\in K:g(x)=0\}$. Then, there exist $N,C\in\bR_{>0}$ such that $\abs{g(x)}^N\leq C\cdot\abs{f(x)}$ for all $x\in K$.

\subsubsection{Corollary}\label{Lojasiewicz_Theorem_general}Let $\mathcal{R}$ be a polynomially bounded o-minimal structure, let
$0<a<1$, and let $f:[0,a]\to\bR$ be a continuous function definable in $\mathcal{R}$. Assume that $f(0)=0$ and $f(\epsilon)>0$ for any $\epsilon\in(0,a]$.
Then, there exist $C>0$ and $m\in\mathbb{N}$, such that $f(\epsilon)\geq
C\cdot\epsilon^m$ for any $\epsilon\in(0,a]$.

\

Proof: This is immediate from Theorem \ref{Lojasiewicz_inequality}. \qed

\subsubsection{Corollary (compact definable sets are globally regular situated)}\label{Lojasiewicz_Theorem_for_definable}Let $\mathcal{R}$
be a polynomially bounded o-minimal structure and let $X,Y\subset \bR^n$ be
compact definable subsets, and assume $X\cap Y\neq\emptyset$. Then, there exist an open neighborhood $V\subset\bR^n$ of $X\cap Y$ and $C,r\in\bR_{>0}$ such
that for all $x\in V$: $$\text{dist}(x,X)+\text{dist}(x,Y)\geq
c\cdot\text{dist}(x,X\cap Y)^r.$$
\

Proof: This follows from Theorem \ref{Lojasiewicz_inequality} -- define on $\{x\in \bR^n:\text{dist}(x,X\cap
Y)\leq2\}$ two functions $f(x):=\text{dist}(x,X)+\text{dist}(x,Y)$, $g(x):=\text{dist}(x,X\cap Y)$  and $V=\{x\in \bR^n:\text{dist}(x,X\cap Y)<1\}$. \qed

\subsection{Exponential structures}A structure $\mathcal{R}$ is called \emph{exponential} if the function $x\mapsto e^x$ (from $\bR$ to $\bR$) is definable in $\mathcal{R}$.

\subsubsection{Theorem (Miller's Dichotomy, \cite{Miller})}\label{Miller_thm}Any o-minimal structure that is not polynomially bounded is exponential.

\subsection{Fr\'echet spaces}\label{prelim-frechet} A Fr\'echet space is
a metrizable, complete, locally convex topological vector space. It can be
shown that the topology of a Fr\'echet space can always be defined by a countable
family of semi-norms. We will use the following results:

\subsubsection{Theorem (closed graph -- cf. \cite[Chapter 17, Corollary 4]{T})}\label{closed_graph_theorem}Let
$E$ and $F$ be Fr\'echet spaces, and $f:E\to F$ a linear map. If the graph
of $f$ inside $E\times F$ is closed, then $f$ is continuous.

\subsubsection{Proposition (cf. \cite[Chapter 10]{T})}\label{closed-sub-of-F-is-F}
A closed subspace of a Fr\'echet space is a Fr\'echet space (with the induced
topology).

\subsubsection{Theorem (Hahn-Banach -- cf. \cite[Chapter 18]{T})}\label{Hahn-Banach}
Let $F$ be a Fr\'echet space, and $K\subset F$ a closed subspace. By Proposition
\ref{closed-sub-of-F-is-F} $K$ is a Fr\'echet space (with the induced topology).
Define $F^*$ (respectively $K^*$) to be the space of continuous linear functionals
on $F$ (on $K$). Then the restriction map $F^*\to K^*$ is onto.

\subsection{Schwartz functions on $\bR^n$}\label{prleim_on_Schwartz_on_Rn}A
Schwartz function on $\bR^n$ is a smooth function $f:\bR^n\to\bR$ such that
for any two multi-indices $\alpha,k\in(\mathbb{N}_0)^n$: $\sup\limits_{x\in\bR^n}\abs{x^k\cdot
f^{(\alpha)}(x)}<\infty$. The space of all Schwartz functions on $\bR^n$
is denoted by $\mathcal{S}(\bR^n)$.

\subsubsection{Proposition (e.g. \cite[Corollary 4.1.2]{AG})}\label{AG-Cor-4.1.2}$\mathcal{S}(\bR^n)$
has a natural structure of a Fr\'echet space, where the topology is given
by the the family of semi-norms indexed by $(\mathbb{N}_0)^n\times(\mathbb{N}_0)^n$:
$$\abs{f}_{\alpha,k}:=\sup\limits_{x\in\bR^n}\abs{x^k\cdot f^{(\alpha)}(x)}.$$

\

We will also use the following result from classic analysis:

\subsection{Proposition (cf. \cite[Corollary 1.4.11]{H})}\label{Hormander's-lemma}Let
$A_0,A_1\subset\bR^n$ be closed subsets. Then, there exists $\phi\in C^\infty(\bR^n\setminus(A_0\cap
A_1))$ such that $\phi|_{A_0\setminus(A_0\cap A_1)}=0$, $\phi|_{A_1\setminus(A_0\cap
A_1)}=1$ and for any multi-index $\alpha\in(\mathbb{N}_0)^n$ there exists
a positive constant $C_\alpha\in\bR$ such that: $$\abs{\phi^{(\alpha)}(x)}\leq
C_\alpha\cdot d(x)^{-|\alpha|}$$ for any $x\in\bR^n\setminus(A_0\cap A_1)$,
where $d(x)=\max\limits_{i=0}^{1}\{\text{dist}(x,A_i)\}$.

\

We end this section with  reminder on the standard stereographic projection:

\subsection{Stereographic projection}\label{stereo_projection}Denote the
unit $n$-dimensional sphere $\{x\in\bR^{n+1}:\abs{x}=1\}$ by $S^n$, and the
standard stereographic (rational bijective) projection from $S^n\setminus\{(1,0,\dots,0)\}$
to $\bR^n$ by $p$. Explicitly, $p$ is defined by $$p((x_0,x_1,\dots,x_n)):=\frac{1}{1-x_0}\cdot(x_1,\dots,x_n),$$
and its inverse is a one point compactification of $\bR^n$, which is given
by $$p^{-1}((x_1,\dots,x_n)):=\frac{1}{(\sum\limits_{j=1}^{n}x_j^2)+1}((\sum\limits_{j=1}^{n}x_j^2)-1,2x_1,2x_2,\dots,2x_n).$$

Note that this (one point) compactification of $\bR^n$ is definable in any o-minimal structure. Thus, we may think
of any definable subset $U\subset\bR^n$ as a definable subset $p^{-1}(U)\subset
S^n$. Since $S^n$ is compact, this point of view will be very useful when we will prove that some functions are Schwartz or tempered on such a subset $U$.

\section{Schwartz functions, tempered functions and tempered distributions on arbitrary open subsets of $\bR^n$}\label{sec-arbitrary}

\subsection{Definition} Let $U\subset\bR^n$ be an open subset, and let $z\in
U$ be some point. We say that a smooth function $f:U\to\bR$ is flat at $z$
if $T_z(f)$ -- its Taylor series at $z$ -- is identically zero. If $Z\subset U$ is any subset,
we say that $f$ is flat at $Z$ if it is flat at $z$, for any $z\in Z$. 

\subsection{Definition}\label{schwartz space for euclidean open definition}
Let $U\subset\bR^n$ be an open subset. Define \emph{the space of Schwartz
functions on $U$}: $$\mathcal{S}(U):=\bigcap\limits_{z\in \bR^n\setminus
U}\bigcap\limits_{\alpha\in(\mathbb{N}_0)^n}\{f\in\mathcal{S}(\bR^n)|f^{(\alpha)}(z)=0\}.$$
By Propositions \ref{closed-sub-of-F-is-F} and \ref{AG-Cor-4.1.2}, $\mathcal{S}(U)$
is a Fr\'echet space. Note that there is a natural bijection between $\mathcal{S}(U)$
and the set $$\{f:U\to\bR|\text{Ext}_U^{\bR^n}(f)\in\mathcal{S}(\bR^n) \text{
and } \text{Ext}_U^{\bR^n}(f)\text{ is flat at }\bR^n\setminus U\}.$$ Thus,
we usually think of Schwartz functions on $U$ as a class of smooth real valued
functions on $U$. In this point of view the topology is given by the family
of semi-norms indexed by $((\mathbb{N}_0)^n)^2$: $$\abs{f}_{\alpha,k}:=\sup\limits_{x\in
U}\abs{x^k\cdot f^{(\alpha)}(x)}.$$ A partial derivative of a Schwartz function
is clearly Schwartz as well.

\subsubsection{Remark}\label{remark_on_extension_by_zero_for_Eucl}If $V\subset
U\subset\bR^n$ are open subsets, then $\mathcal{S}(V)$ is a closed subspace
of $\mathcal{S}(U)$. Also, the extension by zero of any function in $\mathcal{S}(V)$
is a function in $\mathcal{S}(U)$ which is flat at $U\setminus V$, and the
restriction to $V$ of any function in $\mathcal{S}(U)$ which is flat at $U\setminus
V$ is a function in $\mathcal{S}(V)$.

\subsubsection{Proposition}\label{prop_on_equiv_def_of_schwartz} Let $f\in
\mathcal{S}(\bR^n)$. If $f|_{\bR^n\setminus U}\equiv0$, then the following
are equivalent:
\begin{enumerate}
  \item $f\in\mathcal{S}(U)$, i.e. $T_{x_0}(f)\equiv 0$ for any $x_0\in\bR^n\setminus
U$;
  \item for any $m\in\mathbb{N}_0$: $$\sup\limits_{x\in U}\abs{\frac{f(x)}{\text{dist}(x,\bR^n\setminus
U)^m}}<\infty.$$
\end{enumerate}

Proof: Clearly (2) implies (1). Assume (1), and fix some $m\in\mathbb{N}_0$
and some $x\in U$. Take some $x_0\in\bR^n\setminus U$ such that $\text{dist}(x,\bR^n\setminus
U)=\abs{x-x_0}$. Restricting $f$ to the line passing through $x$ and $x_0$,
we use the one dimensional Taylor Theorem with Lagrange's remainder, and
get $$f(x)=T^{m-1}_{x_0}(f)(x)+\frac{f^{(m)}(\xi)}{m!}\cdot\abs{x-x_0}^m$$
for some $\xi$ lying on the line between $x$ and $x_0$, where by $f^{(m)}$
we mean the $m^{\text{th}}$ order one dimensional derivative of $f$ restricted
to this line, i.e. $f^{(m)}=(\sum\limits_{i=1}^{n}a_i\cdot \frac{\partial}{\partial
x_i})^{m}f$, where $(a_1,\dots a_n)$ is the unit vector in the direction
$x-x_0$. Similarly $T^{m-1}_{x_0}(f)$ is the one dimensional Taylor polynomial
of order $m-1$ around $x_0$, of $f$ restricted to this line. By (1) this
is reduced to: $$f(x)=\frac{f^{(m)}(\xi)}{m!}\cdot\abs{x-x_0}^m=\frac{f^{(m)}(\xi)}{m!}\cdot
\text{dist}(x,\bR^n\setminus U)^m.$$
Since $f\in\mathcal{S}(\bR^n)$, for any multi-index $\alpha$ there exists
$C_\alpha>0$ such that $\sup\limits_{y\in\bR^n}\abs{f^{(\alpha)}(y)}<C_\alpha$.
Define $\tilde{C}_m:=\sum\limits_{\abs{\alpha}\leq m}C_\alpha$. As $\abs{a_i}\leq1$
for any $1\leq i\leq n$ one easily sees that $\abs{f^{(m)}(y)}\leq \tilde{C}_m$ for any $y\in\bR^n$, and in particular for $y=\xi$.
Note that this bound is independent of the direction $(a_1,\dots a_n)$. Altogether
$\abs{f(x)}\leq\frac{\tilde{C}_m}{m!}\cdot \text{dist}(x,\bR^n\setminus U)^m$.
\qed

\

The assumption $f\in\mathcal{S}(\bR^n)$ made in Proposition \ref{prop_on_equiv_def_of_schwartz}
is necessary: consider a compactly supported $f:\bR\to\bR$ such that $f|_{\bR^\times}\in
C^\infty(\bR^\times)$ and $f(x)=sin(e^{\frac{1}{x^2}})\cdot e^{-\frac{1}{x^2}}$
for all $0<\abs{x}<1$. Then $f|_{\bR^\times}$ satisfies condition (2), however
$f\not\in C^\infty(\bR)$ and in particular $f|_{\bR^\times}\not\in\mathcal{S}({\bR^\times})$.

\subsection{Definition}\label{def-temp-dist-subsets_of_Rn}
Let $U\subset\bR^n$ be an open subset. Define the \emph{space of tempered
distributions on $U$}, denoted by $\mathcal{S}^*(U)$, as the space of continuous
linear functionals on $\mathcal{S}(U)$.

\subsubsection{Remark on terminology}The term \emph{tempered distributions}
as continuous linear functionals on the space of Schwartz functions is standard in the
context of representation theory (see, for instance, \cite{AG,Fr}). The same
term is sometimes used in the context of analysis, as distributions on subsets (in
the classical sense, i.e. linear functionals on the space of compactly supported
functions) that extend to distributions on the embedding space (see, for
instance, \cite{KS1,KS2}). We will later see that in the polynomially bounded o-minimal case,
and when the embedding space is compact, the two notions coincide. In general
these two notions are different.

\subsection{Lemma}\label{temp-dist-lemma}
Let $V\subset U\subset\bR^n$ be open subsets. Then, $\text{Ext}_V^U:\mathcal{S}(V)\hookrightarrow
\mathcal{S}(U)$ is a closed embedding, and the restriction morphism $\mathcal{S}^*(U)\to\mathcal{S}^*(V)$
is onto.

\

Proof: This follows from Remark \ref{remark_on_extension_by_zero_for_Eucl}
and from the Hahn-Banach Theorem (\ref{Hahn-Banach}). \qed

\subsection{Definition}\label{def_of_temp_func_new}Let $U\subset\bR^n$ be
an open subset. A smooth function
$t:U\to\bR$ is called \emph{tempered} if it satisfies the following two conditions:
\begin{enumerate}
  \item for any $\alpha\in(\mathbb{N}_0)^n$ and any $x_0\in\bR^n\setminus
U$ there exist an open neighborhood $x_0\in V\subset\bR^n$ and $m\in\mathbb{N}_0$
such that: $$\sup\limits_{x\in U\cap V}\abs{t^{(\alpha)}(x)\cdot\text{dist}(x,\bR^n\setminus
U)^m}<\infty;$$
  \item for any $\alpha\in(\mathbb{N}_0)^n$ there exist
$m\in\mathbb{N}_0$ and $r>0$ such that for any $x\in U$ with $\abs{x}>r$:
$$\abs{t^{(\alpha)}(x)}<\abs{x}^m.$$
\end{enumerate} The set of all tempered functions on $U$ is denoted $\mathcal{T}(U)$.

\subsubsection{Remark}Definition \ref{def_of_temp_func_new} is a modification
of the definition of tempered functions appearing in \cite{KS1,KS2}. If $U$
is relatively compact then condition (2) always holds and the two definitions
coincide.
Condition (2) may be thought of as the property of being tempered at infinity.

\

The proofs of the following two lemmas are straightforward:

\subsection{Lemma (pre-sheaf property of tempered functions)}\label{temp_functions_is_a_pre_sheaf}Let
$W\subset U\subset\bR^n$ be open subsets, and let $t\in\mathcal{T}(U)$. Then,
$t|_{W}\in\mathcal{T}(W)$.

\subsection{Lemma}\label{lemma-tempered-is-an-algebra}Let $U\subset\bR^n$
be
an open subset. Then, $\mathcal{T}(U)$ together with point-wise addition
and multiplication is an $\bR$-algebra. Moreover, this algebra is closed
under partial derivatives, and any tempered function that
is bounded from below by some positive number is invertible in this
algebra.

\subsubsection{Remark} A strictly positive tempered function is not necessarily
invertible in $\mathcal{T}(U)$, e.g. $e^{-x^2}\in\mathcal{T}(\bR)$ is strictly
positive, however $e^{x^2}\not\in\mathcal{T}(\bR)$. 

\subsection{Proposition}\label{prop_temp_is_like_module_proprty_on_schwartz}Let
$U\subset\bR^n$ be an open subset and $t\in C^\infty(U)$. Then, $t\in\mathcal{T}(U)$
if and only if for any $s\in\mathcal{S}(U)$, also $t\cdot s\in\mathcal{S}(U)$.

\

Proof: Assume $t\in\mathcal{T}(U)$, and fix $s\in\mathcal{S}(U)$. First let
us show that $\text{Ext}_U^{\bR^n}(t\cdot s)\in C^\infty(\bR^n)$. Let $x_0\in\bR^n\setminus
U$. By condition (1) of Definition \ref{def_of_temp_func_new}, there exist
$C\geq0$, open neighborhood $V\subset\bR^n$ of $x_0$ and $m\in\mathbb{N}_0$
such that for any $x\in U\cap V$ we have $\abs{t(x)}\leq\frac{C}{\text{dist}(x,\bR^n\setminus
U)^m}$. By Proposition \ref{prop_on_equiv_def_of_schwartz} there exists $C'\geq0$
such that for any $x\in U$ we have $\abs{s(x)}\leq C'\cdot\text{dist}(x,\bR^n\setminus
U)^{m+1}$. Thus, for any $x\in U\cap V$ we have: $$\abs{t(x)\cdot s(x)}\leq
\abs{\frac{C}{\text{dist}(x,\bR^n\setminus
U)^m}}\cdot\abs{C'\cdot\text{dist}(x,\bR^n\setminus U)^{m+1}}= C\cdot C'\cdot\text{dist}(x,\bR^n\setminus
U)\leq C\cdot C'\cdot\abs{x-x_0},$$ and so $\text{Ext}_U^{\bR^n}(t\cdot s)$
is continuous at $x_0$. The same proof, together with Leibnitz rule and standard
induction, show that it is fact smooth and flat at $x_0$. The fact that indeed
$\text{Ext}_U^{\bR^n}(t\cdot s)\in\mathcal{S}(\bR^n)$ follows easily from
condition (2) of Definition \ref{def_of_temp_func_new} using Leibnitz rule.
We thus showed $t\cdot s\in\mathcal{S}(U)$.

\

Now assume condition (1) of definition \ref{def_of_temp_func_new} does not
hold. Using Leibnitz rule and standard induction we may assume it does not
hold for $\abs{\alpha}=0$. In this case there exists a point $x_0\in\bR^n\setminus
U$ and a sequence $x_i\in U$ converging to $x_0$, such that for any $i$ we
have $\abs{t(x_i)}\geq\frac{1}{\text{dist}(x_i,\bR^n\setminus U)^i}$. Set
$$B_i:=\{x\in\bR^n:\abs{x-x_i}<\frac{\text{dist}(x_i,\bR^n\setminus U)}{2}\}.$$
As $\bar{B_i}\subset U$, by diluting the sequence we may assume that the
balls $B_i$ are pair-wise disjoint. Define a function $s\in C^\infty(\bR^n)$
by $s|_{B_i}:=\frac{1}{\abs{t(x_i)}}\cdot F(\frac{2\cdot\abs{x-x_i}}{\text{dist}(x_i,\bR^n\setminus
U)})$, where $F:(-1,1)\to\bR$ is defined by $F(x):=e^{\frac{1}{x^2-1}}$,
and $s(x)=0$ outside $\bigcup\limits_{i=1}^{\infty}B_i$. It is easy to verify
that actually $s\in\mathcal{S}(\bR^n)$, and that $\text{supp} (s)\subset
U\cup\{x_0\}$. Also it is easy to see that $s$ is flat at $x_0$, and so $s|_U\in\mathcal{S}(U)$.
However, $t(x_i)\cdot s(x_i)=1$ for any $i$, and $x_i$ converges to $x_0$,
thus $t\cdot (s|_U)\notin\mathcal{S}(U)$. Assuming condition (2) of definition
\ref{def_of_temp_func_new} does not hold, one may construct $s\in\mathcal{S}(U)$
such that $t\cdot s\notin\mathcal{S}(U)$ in a similar way. For instance,
one may apply the mapping $x\mapsto\frac{x}{\abs{x}^2}$ on $U\cap\{\abs{x}>1\}$.
\qed

\subsubsection{Remark}Any tempered function $t\in\mathcal{T}(U)$ defines
a tempered distribution $\xi_t\in\mathcal{S}^*(U)$ by $\xi_t(s):=\int\limits_{x\in
U} (s\cdot t)(x) dx$ for any $s\in\mathcal{S}(U)$. The fact that a smooth
function defines a tempered distribution by integration does not imply that
it is tempered, e.g. $sin(e^x)$ on $\bR$. Not all tempered distributions
arise in such a manner, for instance Dirac's Delta function at some point
$p\in U$, denoted by $\delta_p$, which is defined by $\delta_p(s):=s(p)$
for any $s\in\mathcal{S}(U)$.

\subsection{Lemma}\label{lemma_temp_move_s_to_s}Let $U'\subset U\subset\bR^n$
be open subsets, $t\in\mathcal{T}(U)$, and $s\in\mathcal{S}(U)$. If $\text{supp}(t)\subset
U'$ (where $\text{supp}(t)$ is the Euclidean closure of $t^{-1}(0)$ in $U$, not in
$\bR^n$), then $(t\cdot s)|_{U'}\in\mathcal{S}(U')$.

\

Proof: By Proposition \ref{prop_temp_is_like_module_proprty_on_schwartz}
$t\cdot s\in\mathcal{S}(U)$. Let $x_0\in\bR^n\setminus U'$. We need to show
that $\text{Ext}_U^{\bR^n}(t\cdot s)$ is flat at $x_0$. If $x_0\in\bR^n\setminus
U$ it is obvious. If $x_0\in U\setminus U'$ then, as $t$ is supported on
$U'$, $t$ it is identically zero in a neighborhood of $x_0$, and again $\text{Ext}_U^{\bR^n}(t\cdot
s)$ is flat at $x_0$. \qed

\section{The case of open subsets definable in polynomially bounded o-minimal structures}\label{sec-embedded}

Throughout this section fix $\mathcal{R}$ to be some polynomially bounded o-minimal structure. 

\subsection{Lemma (tempered Urysohn)}\label{tempered_can_separate_disjont_closed}
Let $A_0,A_1\subset U$ be disjoint closed subsets definable in $\mathcal{R}$, where $U\subset\bR^n$
is some open definable subset. Then, there exists a non-negative $t\in\mathcal{T}(U)$
such that $t|_{A_0}=0$ and $t|_{A_1}=1$.

\

Proof: Denote by $\tilde{U}$ (respectively $\tilde A_i$) the pre-image of
$U$ (resp. $A_i$) in $S^n$ under the standard stereographic projection (see
\ref{stereo_projection}), and denote by $\bar{A}_i$ the closure of $\tilde
A_i$ inside $S^n$. Then, $\bar{A}_0$ and $\bar{A}_1$ are closed definable
subsets of $\bR^{n+1}$. Note that as $\tilde A_0$ and $\tilde A_1$ are disjoint,
$\bar{A}_0\cap \bar{A}_1\cap\tilde{U}=\emptyset$. Before we proceed replace
$\bar{A}_0$ by $\bar{A}_0\cup(S^n\setminus\tilde{U})$ -- this does not change
neither the fact that $\bar{A}_0$ is a closed definable subset of $\bR^{n+1}$,
the fact that $\bar{A}_0\cap\tilde{U}=\tilde{A}_0$ nor the fact that $\bar{A}_0\cap
\bar{A}_1\cap\tilde{U}=\emptyset$.

\

By Proposition \ref{Hormander's-lemma} there exists $\phi\in C^\infty(\bR^{n+1}\setminus(\bar
A_0\cap \bar A_1))$ such that $\phi|_{\bar A_0\setminus(\bar A_0\cap \bar
A_1)}=0$, $\phi|_{\bar A_1\setminus(\bar A_0\cap \bar A_1)}=1$ and for any
$\alpha\in(\mathbb{N}_0)^{n+1}$ there exists a positive constant $C_\alpha\in\bR$
such that: $$(1)\text{ }\text{ }\text{ }\text{ }\abs{\phi^{(\alpha)}(x)}\leq
C_\alpha\cdot d(x)^{-|\alpha|},$$ for any $x\in\bR^{n+1}\setminus(\bar A_0\cap
\bar A_1)$, where $d(x)=\max\limits_{i=0}^{1}\{\text{dist}(x,\bar A_i)\}$.
Note that $\bar{A}_0\cap \bar{A}_1\cap\tilde{U}=\emptyset$ implies that $\tilde{U}\subset\bR^{n+1}\setminus(\bar
A_0\cap \bar A_1)$, and so $\phi$ is defined on $\tilde{U}$. We may of course
assume $\phi$ vanishes whenever $\abs{x}>3$, otherwise we replace it by a
different smooth function without changing its values in any $\abs{x}<2$.

\

By Corollary \ref{Lojasiewicz_Theorem_for_definable} there exist $c,r>0$
and an open neighborhood $V\subset\bR^{n+1}$ of $\bar{A}_0\cap\bar{A}_1$,
such that for any $x\in V$: $$(2)\text{ }\text{ }\text{ }\text{ }\text{ }\text{
}\text{ }\text{ }d(x)\geq\frac{1}{2}(\text{dist}(x,\bar{A}_0)+\text{dist}(x,\bar{A}_1))\geq\frac{c}{2}\cdot\text{dist}(\bar{A}_0\cap\bar{A}_1)^r.$$
Combining (1) and (2) we get that for $x\in V$: $$\abs{\phi^{(\alpha)}(x)}\leq
C_\alpha\cdot d(x)^{-|\alpha|}\leq\frac{C_\alpha\cdot c}{2}\cdot\text{dist}(x,\bar
A_0\cap \bar A_1)^{-\abs{\alpha}\cdot r}\leq\frac{C_\alpha\cdot c}{2}\cdot\text{dist}(x,\bar
A_0\cap \bar A_1)^{-m},$$ for any $m\in\mathbb{N}$ satisfying $m>\abs{\alpha}\cdot
r$. We conclude that $\phi\in\mathcal{T}(\bR^{n+1}\setminus(\bar A_0\cap
\bar A_1))$.

\

Denote by $\check{\phi}$ the push-forward of $\phi|_{\tilde U}$ under the
stereographic projection, i.e. $\check{\phi}(x)=\phi(p^{-1}(x))$. As the
distances on $S^n$ and the distances on $\bR^n$ of corresponding points under
the stereographic projection are related by rational functions, and as the
Jacobian of the stereographic projection is also rational, we find that $\check{\phi}\in\mathcal{T}(U)$,
where $\check{\phi}|_{A_0}=0$ and $\check{\phi}|_{A_1}=1$. In order to get
a non-negative $t\in\mathcal{T}(U)$ with the same properties, apply Lemma
\ref{lemma-tempered-is-an-algebra} and define $t(x)=\check{\phi}(x)^2$. \qed

\subsection{Proposition (tempered partition of unity on definable open subsets of $\bR^n$)}\label{tempered-partition-of-unity}Let $U\subset\bR^n$
be an open definable subset and $U=\bigcup\limits_{i=1}^{m}U_i$ be a finite
open definable cover. Then, there exist non-negative tempered functions
$\varphi_1,\dots,\varphi_m\in\mathcal{T}(U)$ such that the following hold:\begin{enumerate}
                                                       \item $\text{supp}(\varphi_i)\subset
U_i$ for any $i=1,\dots,m$ (where $\text{supp}(\varphi_i)$ is the closure
of $\varphi_i^{-1}(\bR\setminus\{0\})$ inside $U$, not inside $\mathbb{R}^n$);
                                                       \item $\sum\limits_{i=1}^{m}\varphi_i(x)=1$
for any $x\in U$.
                                                     \end{enumerate}

\

Proof: Our proof is based on the idea of proving continuous partition of
unity for arbitrary normal topological spaces exists, as given in \cite[Thorem
36.1]{Mun}, with some modifications. For the sake of completeness we repeat
all the details:

\

Normality (Lemma \ref{normality_of_definable_lemma}) implies that for
the open cover $U_1,U_2,\dots,U_m$ we can find an open definable cover
$V_1,V_2,\dots,V_m$ of $U$, such that $\overline{V_i}\subset U_i$ for any
$i=1,2,\dots,m$ (this is sometimes called "the shrinking lemma", and its
proof -- which is made by standard induction -- is given in the beginning
of the proof of \cite[Thorem 36.1]{Mun}). Applying this procedure again we
can find an open definable cover $W_1,W_2,\dots,W_m$ of $U$ such that
$\overline{W_i}\subset V_i$ for any $i=1,2,\dots,m$.

\

For any $i=1,2,\dots,m$ we apply Lemma \ref{tempered_can_separate_disjont_closed}
to the disjoint closed sets $\overline{W_i},U\setminus V_i\subset U$, and
get a non-negative function $t_i\in\mathcal{T}(U)$ such that $t_i|_{\overline{W_i}}=1$
and $t_i|_{U\setminus V_i}=0$. Since $t_i^{-1}(\mathbb{R}\setminus\{0\})\subset
V_i$, we have $\text{supp}(t_i)\subset\overline{V_i}\subset U_i$. Note that
by Lemma \ref{lemma-tempered-is-an-algebra} the function $\sum\limits_{i=1}^{m}t_i$
is a tempered function on $U$, and as it is bounded from below by $1$, its
inverse is a tempered function as well. Finally define $\varphi_i:=\frac{t_i}{\sum\limits_{i=1}^{m}t_i}$.
Clearly these $\varphi_i$'s satisfy (1) and (2). \qed

\subsection{Proposition (sheaf property of tempered functions on definable open subsets
of $\bR^n$)}\label{tempered-is-a-sheaf-definable-subset-of-Rn}
Let $U\subset\bR^n$ be an open definable subset and $U=\bigcup\limits_{i=1}^{m}U_i$
be a finite open definable cover. Then, the following sequence is exact:

\begin{center}
\leavevmode
\xymatrix{
0 \ar[rr] & & \mathcal{T}(U) \ar[rr]^{\text{Res}_1} & & \bigoplus\limits_{i=1}^{m}\mathcal{T}(U_i)
\ar[rr]^(.4){\text{Res}_2} & & \bigoplus\limits_{i=1}^{m-1}\bigoplus\limits_{j=i+1}^{m}\mathcal{T}(U_i\cap
U_j),}
\end{center}
where the $k$-th coordinate of $\text{Res}_1(t)$ is $\text{Res}_{U_k}^{U}(t)$
and for $1\leq k< l\leq m$ the $(k,l)$-coordinate of $\text{Res}_2(\bigoplus\limits_{i=1}^{m}t_i)$
is $\text{Res}_{U_k\cap U_l}^{U_k}(t_k)-\text{Res}_{U_k\cap U_l}^{U_l}(t_l)$.

\

Proof: It follows from Lemma \ref{temp_functions_is_a_pre_sheaf} that the
above maps are well defined. It is enough to prove the following: let $t_i\in\mathcal{T}(U_i)$
be such that $t_i|_{U_i\cap U_j}=t_j|_{U_i\cap U_j}$ for any $1\leq i\leq j\leq
m$, and define $t:U\to \bR$ by $t|_{U_i}=t_i$ for any $i=1,2,\dots,m$. Then
$t\in\mathcal{T}(U)$. By Proposition \ref{prop_temp_is_like_module_proprty_on_schwartz}
it is enough to show that $t\cdot s\in\mathcal{S}(U)$ for any $s\in\mathcal{S}(U)$.
Indeed, let $s\in\mathcal{S}(U)$. By Proposition \ref{tempered-partition-of-unity}
there exist non-negative tempered functions $\varphi_1,\dots,\varphi_m\in\mathcal{T}(U)$
such that $\text{supp}(\varphi_i)\subset U_i$ for any $i=1,\dots,m$, and
$\sum\limits_{i=1}^{m}\varphi_i(x)=1$ for any $x\in U$. Lemma \ref{lemma_temp_move_s_to_s}
implies that $\varphi_i|_{U_i}\cdot s|_{U_i}\in\mathcal{S}(U_i)$, and so,
again by Proposition \ref{prop_temp_is_like_module_proprty_on_schwartz},
$t_i\cdot \varphi_i|_{U_i}\cdot s|_{U_i}\in\mathcal{S}(U_i)$. Recalling that
the extension by zero of a Schwartz function from the open subset $U_i$ is
a Schwartz function on $U$, and that finite sums of Schwartz functions are
Schwartz functions, we have $t\cdot s=\sum\limits_{i=1}^{m}\text{Ext}_{U_i}^{U}(t_i\cdot
\varphi_i|_{U_i}\cdot s|_{U_i})\in\mathcal{S}(U)$. \qed

\

The following Proposition \ref{tempered-distributions-is-a-sheaf-definable-subset-of-Rn}
follows from Remark \ref{remark_on_extension_by_zero_for_Eucl}, Lemma \ref{temp-dist-lemma},
Lemma \ref{lemma_temp_move_s_to_s}, Proposition \ref{tempered-partition-of-unity}
and Proposition \ref{tempered-is-a-sheaf-definable-subset-of-Rn}. The
following Proposition \ref{schwartz_functions-is-a-cosheaf-definable-subset-of-Rn}
follows from Lemma \ref{temp-dist-lemma} and Proposition \ref{tempered-is-a-sheaf-definable-subset-of-Rn}.
The proofs, which we omit, are quite technical, and are identical to the
proofs of \cite[Proposition 4.4 and Proposition 4.5]{ES} (respectively),
with straightforward adjustments.

\subsection{Proposition (sheaf property of tempered distributions on definable open subsets
of $\bR^n$)}\label{tempered-distributions-is-a-sheaf-definable-subset-of-Rn}
Let $U\subset\bR^n$ be an open definable subset and $U=\bigcup\limits_{i=1}^{m}U_i$
be a finite open definable cover. Then, the following sequence is exact:

\begin{center}
\leavevmode
\xymatrix{
0 \ar[rr] & & \mathcal{S}^*(U) \ar[rr]^{\text{Res}_1} & & \bigoplus\limits_{i=1}^{m}\mathcal{S}^*(U_i)
\ar[rr]^(.4){\text{Res}_2} & & \bigoplus\limits_{i=1}^{m-1}\bigoplus\limits_{j=i+1}^{m}\mathcal{S}^*(U_i\cap
U_j),}
\end{center}
where the $k$-th coordinate of $\text{Res}_1(\xi)$ is $\xi|_{\mathcal{S}(U_k)}$
and for $1\leq k< l\leq m$ the $(k,l)$-coordinate of $\text{Res}_2(\bigoplus\limits_{i=1}^{m}\xi_i)$
is $\xi_k|_{\mathcal{S}(U_k\cap U_l)}-\xi_l|_{\mathcal{S}(U_k\cap U_l)}$.

\subsection{Proposition (cosheaf property of Schwartz functions on definable
open subsets
of $\bR^n$)}\label{schwartz_functions-is-a-cosheaf-definable-subset-of-Rn}
Let $U\subset\bR^n$ be an open definable subset and $U=\bigcup\limits_{i=1}^{m}U_i$
be a finite open definable cover. Then, the following sequence is exact:

\begin{center}
\leavevmode
\xymatrix{
\bigoplus\limits_{i=1}^{m-1}\bigoplus\limits_{j=i+1}^{m}\mathcal{S}(U_i\cap
U_j) \ar[rr]^(.6){\text{Ext}_1} & & \bigoplus\limits_{i=1}^{m}\mathcal{S}(U_i)
\ar[rr]^{\text{Ext}_2} & & \mathcal{S}(U) \ar[rr] & & 0,}
\end{center}
where the $k$-th coordinate of $\text{Ext}_1(\bigoplus\limits_{i=1}^{m-1}\bigoplus\limits_{j=i+1}^{m}s_{i,j})$
is $\sum\limits_{i=1}^{k-1} \text{Ext}_{U_k\cap U_i}^{U_k}(s_{i,k})-\sum\limits_{i=k+1}^{m}
\text{Ext}_{U_k\cap U_i}^{U_k}(s_{k,i})$, and $\text{Ext}_2(\bigoplus\limits_{i=1}^{m}s_i):=\sum\limits_{i=1}^{m}\text{Ext}^{U}_{U_i}(s_i)$.

\

We end this section by proving that a definable diffeomorphism (in a polynomially bounded o-minimal structure) between
open definable subsets of $\bR^n$ induces an isomorphism of the corresponding
Schwartz spaces (Corollary \ref{definable-diffeo-induce-iso-of-schwartz}).
We will need some lemmas before.

\subsection{Lemma}\label{definable-smooth-function-is-tempered}Let $U\subset\bR^n$
be an open definable subset and $f:U\to\bR$ be a definable smooth function.
Then, $f\in\mathcal{T}(U)$.

\

Proof: To avoid problems in infinity, we again consider $U$ as an open subset
of $S^n$ (see the proof of Lemma \ref{tempered_can_separate_disjont_closed}
for more detail). If $f$ is identically zero there is nothing to prove. Assume
$f$ is not identically zero. As partial derivative of a definable function
is definable, it is enough to prove that there exist $\epsilon',C>0$ and
$m\in\mathbb{N}_0$ such that $\abs{f(x)}\leq C\cdot\text{dist}(x,\partial
U)^{-m}$ for any $x\in U$ satisfying $\text{dist}(x,\partial U)<\epsilon'$.

\

Define a definable function $$F(\epsilon):=\max\limits_{\{x:\text{dist}(x,\partial
U)\geq\epsilon\}}\abs{f(x)},$$ for all $\epsilon\in(0,\epsilon']$, where
$0<\epsilon'<1$ is sufficiently small such that $F(\epsilon')>0$. Note that
$F$ is continuous and monotonic decreasing. Define a continuous definable
function $G(\epsilon):=\frac{1}{F(\epsilon)}$ for all $\epsilon\in(0,\epsilon']$.
Note that $G$ is monotonic increasing, and so $G(0):=\lim\limits_{\epsilon\to0}G(\epsilon)$
exists, and non-negative. If $G(0)>0$ then $f$ is bounded on $U$, and we
are done. Other-wise $G(0)=0$, and $G(\epsilon)>0$ for any $\epsilon>0$.
Applying Corollary \ref{Lojasiewicz_Theorem_general}, there exist $C>0$ and
$m\in\mathbb{N}$, such that $G(\epsilon)\geq C\cdot\epsilon^m$ for any $\epsilon\in(0,\epsilon']$.
We conclude that for any $\epsilon\in(0,\epsilon']$: $$\max\limits_{\{x:\text{dist}(x,\partial
U)\geq\epsilon\}}\abs{f(x)}\leq C^{-1}\cdot\epsilon^{-m}.$$
Finally, for $x_0\in U$ such that $\epsilon_0:=\text{dist}(x_0,\partial U)<\epsilon'$,
we have: $$\abs{f(x_0)}\leq\max\limits_{\{x:\text{dist}(x,\partial U)\geq\epsilon_0\}}\abs{f(x)}\leq
C^{-1}\cdot\epsilon_0^{-m}=C^{-1}\cdot\text{dist}(x_0,\partial U)^{-m}.$$
\qed

\subsection{Lemma}\label{definable-diffeo-pull-schwartz-to-schwartz}Let
$U,V\subset\bR^n$ be open definable subsets, and $\nu:U\to V$ a definable
diffeomorphism. Then, $\nu^*(\mathcal{S}(V))\subset\mathcal{S}(U)$.

\

Proof: We start with the one dimensional case: let $s\in\mathcal{S}(V)$,
we need to show that $\nu^*(s)\in\mathcal{S}(U)$. Let $x_0\in\partial U$.
We first wish to prove that for any $l\in\mathbb{N}_0$ we have $\lim\limits_{x\to
x_0}(\nu^*(s))^{(l)}(x)=0$. We only show the cases $l=0,1$ -- the rest is
left to the reader. Assume toward a contradiction there exist $\epsilon>0$
and $x_1,x_2,x_3,\dots\in U$ such that $\lim\limits_{n\to\infty}x_n=x_0$
and for any $n\in\mathbb{N}$ we have $\abs{s(\nu(x_n))}>\epsilon$. Denote
$y_n:=\nu(x_n)\in V$, and denote by $y_0\in\partial V$ some accumulation
point of the sequence $\{y_n\}_{n=1}^{\infty}$ if such exists, and $y_0=\pm\infty$
otherwise. In any case by diluting the sequence $\{x_n\}_{n=1}^{\infty}$
we may assume $\lim\limits_{n\to\infty}y_n=y_0$. Then, $\epsilon\leq\lim\limits_{x\to
x_0}\abs{s(\nu(x_n))}=\lim\limits_{y\to y_0}\abs{s(y_n)}=0$, as $y_0\not\in
V$ and $s\in\mathcal{S}(V)$. This is a contradiction, and so the case $l=0$
is done.

\

Now assume toward a contradiction there exist $\epsilon>0$ and $x_1,x_2,x_3,\dots\in
U$ such that $\lim\limits_{n\to\infty}x_n=x_0$ and for any $n\in\mathbb{N}$
we have $\abs{\nu^*(s)'(x_n)}>\epsilon$. Denote $y_n:=\nu(x_n)\in V$, and
denote by $y_0\in\partial V$ some accumulation point of the sequence $\{y_n\}_{n=1}^{\infty}$
if such point exists, and $y_0=\pm\infty$ otherwise. As before, by moving
to a subsequence we may assume $\lim\limits_{n\to\infty}y_n=y_0$. Then, $$\epsilon\leq\lim\limits_{x\to
x_0}\abs{(\nu^*(s))'(x_n)}=\lim\limits_{x\to x_0}\abs{s'(\nu(x_n))\cdot \nu'(x_n)}=$$$$\lim\limits_{y\to
y_0}\abs{s'(y_n)\cdot \nu'(\nu^{-1}(y_n))}=\lim\limits_{y\to y_0}\abs{s'(y_n)\cdot
(\nu^{-1})^*(\nu')(y_n)}=0,$$ as $y_0\not\in V$, $s'\in\mathcal{S}(V)$, $(\nu^{-1})^*(\nu')\in\mathcal{T}(V)$
(as it is definable and by Lemma \ref{definable-smooth-function-is-tempered}),
and so $s'\cdot(\nu^{-1})^*(\nu')\in\mathcal{S}(V)$ (by Proposition \ref{prop_temp_is_like_module_proprty_on_schwartz}).
This is a contradiction, and the case  $l=1$ is done.

\

Thus we showed that $\text{Ext}_{U}^{\bR}(\nu^*(s))\in C^\infty(\bR)$ and
that it is flat outside $U$. In order to show that it is a Schwartz function,
we are left to check it also behaves well at infinity, which is essentially
checking one more point on the circle. This can be done in the same way as
we did for usual points in $\bR$, and is left for the reader to verify.

\

The general $n$-dimensional case follows in the same way, but now one has
to apply a general $n$-dimensional algebraic differential operator $\mathcal{L}$
on $\nu^*(s)$, instead of a one-dimensional derivative. One has to use the
following fact: let $U\subset\bR^n$, $V\subset\bR^m$ be two open definable
subsets, and let $\nu:U\to V$ be some map, given by $\nu(x)=(\nu_1(x),\dots,\nu_m(x))$,
for some functions $\nu_1,\dots,\nu_m:U\to\bR$. Then, $\nu$ is a definable
map if and only if the functions $\nu_1,\dots,\nu_m$ are all definable.
All the other details may be verified straightforward and are left to the
reader. \qed

\subsection{Lemma}\label{enough-to-pull-Eucl}Let $U,V\subset\bR^n$ be open
subsets (not necessarily definable), and $\nu:U\to V$ some map (not necessarily definable). If $\nu^*(\mathcal{S}(V))\subset\mathcal{S}(U)$,
then $\nu^*:\mathcal{S}(V)\to\mathcal{S}(U)$ is continuous.

\

Proof: $$\mathcal{S}(V)\times\mathcal{S}(U)\supset \text{graph}(\nu^*)=$$$$\{(f,\nu^*f)|f\in\mathcal{S}(V)\}=\{(f,g)|g(x)=
f(\nu(x))\text{ }\forall x\in U\}=$$$$\bigcap\limits_{x\in U}\{(f,g)|g(x)=
f(\nu(x))\}.$$ For any $x_0\in U$ denote $$\mathcal{S}(V)\times\mathcal{S}(U)\supset
A_{x_0}:=\{(f,g)|g({x_0})= f(\nu({x_0}))\}.$$ Let $(f,g)\in\mathcal{S}(V)\times\mathcal{S}(U)\setminus
A_{x_0}$, i.e. $\epsilon:=\abs{g({x_0})-f(\nu({x_0}))}>0$. For any $\tilde{f}\in\mathcal{S}(V)$
such that $\abs{f-\tilde{f}}_{(0,0)}:=\sup\limits_{y\in V}\abs{f(y)-\tilde{f}(y)}<\frac{\epsilon}{4}$
and any $\tilde{g}\in\mathcal{S}(U)$ such that $\abs{g-\tilde{g}}_{(0,0)}:=\sup\limits_{x\in
U}\abs{g(x)-\tilde{g}(x)}<\frac{\epsilon}{4}$, we also have $(\tilde{f},\tilde{g})\in\mathcal{S}(V)\times\mathcal{S}(U)\setminus
A_{x_0}$. Thus $A_{x_0}$ is closed in $\mathcal{S}(V)\times\mathcal{S}(U)$,
and so $\text{graph}(\nu^*)=\bigcap\limits_{x_0\in U}A_{x_0}$ is closed in
$\mathcal{S}(V)\times\mathcal{S}(U)$. By Theorem \ref{closed_graph_theorem}
$\nu^*$ is continuous. \qed

\subsection{Corollary}\label{definable-diffeo-induce-iso-of-schwartz}Let
$U,V\subset\bR^n$ be open definable subsets, and $\nu:U\to V$ a definable
diffeomorphism. Then, $\nu^*|_{\mathcal{S}(V)}$ is an isomorphism of the
Fr\'echet spaces $\mathcal{S}(V)$ and $\mathcal{S}(U)$.

\

Proof: It follows directly from Lemma \ref{definable-diffeo-pull-schwartz-to-schwartz}
and Lemma \ref{enough-to-pull-Eucl}. \qed

\section{The affine polynomially bounded case}\label{sec-affine}

Throughout this section fix $\mathcal{R}$ to be some polynomially bounded
o-minimal structure. 

\subsection{Definition}Let $M$ be a smooth definable manifold and let $f$ be a real valued function defined on an open definable neighborhood
of a point $z\in M$. $f$
is called \emph{flat at $z\in M$} if there exists a chart $(U,g)$ with $z\in
U$ such that (or equivalently if for any such a chart) $g_*(f):=f\circ g^{-1}$
is flat
at $g(z)$. It is called flat at $Z\subset M$ if it is flat at any $z\in
Z$. Note that this definition makes sense for general smooth manifolds (not necessarily definable), however we do not need this generality in this paper.

\subsection{Definition}An $n$-dimensional definable manifold $M$ is called
\emph{affine} if there exists (in one of the atlases in the appropriate equivalent
class of atlases) a chart of the form $(M,g)$, with $g(M)\subset\bR^n$ is open.

\subsection{Definitions and easy observations}\label{definitions-for-affine}Let
$M$ be an affine manifold definable in $\mathcal{R}$.

\begin{enumerate}

\item A function $f:M\to\bR$ is called \emph{a Schwartz function on $M$}
if $g_*f:=f\circ
g^{-1}\in\mathcal{S}(g(M))$, where $(M,g)$ is a chart on $M$ with $g(M)\subset\bR^n$ is
open. Denote the
space of all Schwartz
functions on $M$ by $\mathcal{S}(M)$, and define a topology on $\mathcal{S}(M)$
by declaring a subset $U\subset\mathcal{S}(M)$ to be open if $g_*(U)\subset\mathcal{S}(g(M))$
is an open subset. If $(M,g')$ is another chart in the same equivalence class
of atlases with $g'(M)\subset\bR^n$ is
open, then Corollary \ref{definable-diffeo-induce-iso-of-schwartz}
implies
that $g'\circ g^{-1}$ is an isomorphism of $\mathcal{S}(g(M))$
and $\mathcal{S}(g'(M))$, and so $\mathcal{S}(M)$ is a well defined
Fr\'echet space.
Note that by Lemma \ref{temp-dist-lemma} if $V\subset
U\subset M$ are open definable subsets, then $\mathcal{S}(V)$ is a closed
subspace
of $\mathcal{S}(U)$: explicitly (by Remark \ref{remark_on_extension_by_zero_for_Eucl}) $\text{Ext}_V^U:\mathcal{S}(V)\hookrightarrow
\mathcal{S}(U)$ is a closed embedding, whose image consists of all functions
in $\mathcal{S}(U)$ which are flat at $U\setminus V$.

\item The \emph{space of tempered
distributions on $M$}, denoted by $\mathcal{S}^*(M)$, is the space of continuous
linear functionals on $\mathcal{S}(M)$. By Lemma \ref{temp-dist-lemma} if
$V\subset U\subset M$ are open definable subsets, then the restriction morphism
$\mathcal{S}^*(U)\to\mathcal{S}^*(V)$
is onto.

\item A function
$f:M\to\bR$ is called \emph{tempered} if for any $s\in\mathcal{S}(M)$ also
$f\cdot
s\in\mathcal{S}(M)$. The set of all tempered functions on $M$ is denoted
$\mathcal{T}(M)$. Note that Proposition \ref{prop_temp_is_like_module_proprty_on_schwartz}
implies that $f\in\mathcal{T}(M)$ if and only if $g_*f:=f\circ
g^{-1}\in\mathcal{T}(g(M))$, where $(M,g)$ is a chart on $M$ with $g(M)\subset\bR^n$ is
open. Also note that
(by  Lemma \ref{lemma-tempered-is-an-algebra}) $\mathcal{T}(M)$ together with point-wise addition and multiplication is an $\bR$-algebra
and any tempered function that
is bounded from below by some positive number is invertible in this
algebra.

\end{enumerate}

The following results follow immediately from the corresponding results on
open definable subsets of $\bR^n$:

\subsection{Lemma (cf. Lemma \ref{lemma_temp_move_s_to_s})}\label{lemma_temp_move_s_to_s_affine}Let
$M$ be an affine definable manifold, $U\subset M$ an open definable
subset, $t\in\mathcal{T}(M)$, and $s\in\mathcal{S}(M)$. If $\text{supp}(t)\subset
U$, then $(t\cdot s)|_{U}\in\mathcal{S}(U)$.

\subsection{Proposition (tempered partition of unity on affine definable manifolds -- cf. Proposition \ref{tempered-partition-of-unity})}\label{tempered-partition-of-unity-affine}Let
$M$ be an affine definable manifold and $M=\bigcup\limits_{i=1}^{m}U_i$
be a finite open definable cover. Then, there exist non-negative tempered
functions $\varphi_1,\dots,\varphi_m\in\mathcal{T}(M)$ such that the following
hold:\begin{enumerate}
                                                       \item $\text{supp}(\varphi_i)\subset
U_i$ for any $i=1,\dots,m$;
                                                       \item $\sum\limits_{i=1}^{m}\varphi_i(x)=1$
for any $x\in M$.
                                                     \end{enumerate}

\subsection{Proposition (sheaf property of tempered functions on affine definable
manifolds -- cf. Proposition \ref{tempered-is-a-sheaf-definable-subset-of-Rn})}\label{tempered-is-a-sheaf-affine}
Let $M$ be an affine definable manifold and $M=\bigcup\limits_{i=1}^{m}U_i$
be a finite open definable cover. Then, the following sequence is exact:

\begin{center}
\leavevmode
\xymatrix{
0 \ar[rr] & & \mathcal{T}(M) \ar[rr]^{\text{Res}_1} & & \bigoplus\limits_{i=1}^{m}\mathcal{T}(U_i)
\ar[rr]^(.4){\text{Res}_2} & & \bigoplus\limits_{i=1}^{m-1}\bigoplus\limits_{j=i+1}^{m}\mathcal{T}(U_i\cap
U_j),}
\end{center}
where the $k$-th coordinate of $\text{Res}_1(t)$ is $\text{Res}_{U_k}^{M}(t)$
and for $1\leq k< l\leq m$ the $(k,l)$-coordinate of $\text{Res}_2(\bigoplus\limits_{i=1}^{m}t_i)$
is $\text{Res}_{U_k\cap U_l}^{U_k}(t_k)-\text{Res}_{U_{k}\cap U_l}^{ U_l}(t_l)$.

\subsection{Proposition (sheaf property of tempered distributions on affine definable
manifolds -- cf. Proposition \ref{tempered-distributions-is-a-sheaf-definable-subset-of-Rn})}\label{tempered-distributions-is-a-sheaf-affine}
Let $M$ be an affine definable manifold and $M=\bigcup\limits_{i=1}^{m}U_i$
be a finite open definable cover. Then, the following sequence is exact:

\begin{center}
\leavevmode
\xymatrix{
0 \ar[rr] & & \mathcal{S}^*(M) \ar[rr]^{\text{Res}_1} & & \bigoplus\limits_{i=1}^{m}\mathcal{S}^*(U_i)
\ar[rr]^(.4){\text{Res}_2} & & \bigoplus\limits_{i=1}^{m-1}\bigoplus\limits_{j=i+1}^{m}\mathcal{S}^*(U_i\cap
U_j),}
\end{center}
where the $k$-th coordinate of $\text{Res}_1(\xi)$ is $\xi|_{\mathcal{S}(U_k)}$
and for $1\leq k< l\leq m$ the $(k,l)$-coordinate of $\text{Res}_2(\bigoplus\limits_{i=1}^{m}\xi_i)$
is $\xi_k|_{\mathcal{S}(U_k\cap U_l)}-\xi_l|_{\mathcal{S}(U_k\cap U_l)}$.

\subsection{Proposition (cosheaf property of Schwartz functions on affine
definable
manifolds -- cf. Proposition \ref{schwartz_functions-is-a-cosheaf-definable-subset-of-Rn})}\label{schwartz_functions-is-a-cosheaf-affine}
Let $M$ be an affine definable manifold and $M=\bigcup\limits_{i=1}^{m}U_i$
be a finite open definable cover. Then, the following sequence is exact:

\begin{center}
\leavevmode
\xymatrix{
\bigoplus\limits_{i=1}^{m-1}\bigoplus\limits_{j=i+1}^{m}\mathcal{S}(U_i\cap
U_j) \ar[rr]^(.6){\text{Ext}_1} & & \bigoplus\limits_{i=1}^{m}\mathcal{S}(U_i)
\ar[rr]^{\text{Ext}_2} & & \mathcal{S}(M) \ar[rr] & & 0,}
\end{center}
where the $k$-th coordinate of $\text{Ext}_1(\bigoplus\limits_{i=1}^{m-1}\bigoplus\limits_{j=i+1}^{m}s_{i,j})$
is $\sum\limits_{i=1}^{k-1} \text{Ext}_{U_k\cap U_i}^{U_k}(s_{i,k})-\sum\limits_{i=k+1}^{m}
\text{Ext}_{U_k\cap U_i}^{U_k}(s_{k,i})$, and $\text{Ext}_2(\bigoplus\limits_{i=1}^{m}s_i):=\sum\limits_{i=1}^{m}\text{Ext}^{M}_{U_i}(s_i)$.

\section{The general polynomially bounded case}\label{sec-general}

Throughout this section fix
$\mathcal{R}$ to be some polynomially bounded o-minimal structure.

\

{\bf Some more terminology.} If $M$
is a definable
manifold and $(U,g)$ is a definable chart on $M$ with $g(U)\subset\bR^n$ is
open, then
$U$ has a natural structure of an affine definable manifold induced from $M$.
In that case we call $U\subset M$ an \emph{affine definable patch}. Accordingly,
if $U_1,U_2,\dots,U_m$
are affine definable patches such that $M=\bigcup\limits_{i=1}^{m}U_i$, then we
call $U_1,U_2,\dots,U_m$ an \emph{affine definable open cover of $M$}. Finally,
for any definable manifold $M$, we denote the space of all real valued
functions on $M$ by $\text{Func}(M,\bR)$.

\subsection{Lemma}\label{lemma-schwartz-on-sub-analytic-is-well-defined}Let
$M$ be a definable manifold, and let $M=\bigcup\limits_{i=1}^mU_i=\bigcup\limits_{i=m+1}^lU_i$
be two affine definable open covers. Using extensions by zero, there are natural maps $\phi_1:\bigoplus\limits_{i=1}^m\text{Func}(U_i,\bR)\to
\text{Func}(M,\bR)$ and $\phi_2:\bigoplus\limits_{i=m+1}^l\text{Func}(U_i,\bR)\to
\text{Func}(M,\bR)$. Then, $\phi_1(\bigoplus\limits_{i=1}^m\mathcal{S}(U_i))\cong(\bigoplus\limits_{i=1}^m\mathcal{S}(U_i))/Ker(\phi_1|_{\bigoplus\limits_{i=1}^m\mathcal{S}(U_i)})$
has a natural structure of a Fr\'echet space, and moreover there is an isomorphism
of Fr\'echet spaces $\phi_1(\bigoplus\limits_{i=1}^m\mathcal{S}(U_i))\cong\phi_2(\bigoplus\limits_{i=m+1}^l\mathcal{S}(U_i))$.

\

Proof: This follows from \ref{definitions-for-affine}(1) and Proposition \ref{schwartz_functions-is-a-cosheaf-affine}. The detailed
proof, which we omit, is identical to the proof of \cite[Lemma 5.1]{ES},
with straightforward adjustments. \qed

\

\subsection{Definition}\label{def_schwartz_fun_on_definable}Let $M$ be
a definable manifold, let $M=\bigcup\limits_{i=1}^mU_i$ be some affine
definable open cover, and consider the natural map $\phi:\bigoplus\limits_{i=1}^m\text{Func}(U_i,\bR)\to
\text{Func}(M,\bR)$. Define \emph{the space of Schwartz functions on $M$}
by $\mathcal{S}(M):=(\bigoplus\limits_{i=1}^m\mathcal{S}(U_i))/\text{Ker}(\phi|_{\bigoplus\limits_{i=1}^m\mathcal{S}(U_i)})$,
with the natural quotient topology. By Lemma \ref{lemma-schwartz-on-sub-analytic-is-well-defined}
this definition is independent of the cover chosen, and $\mathcal{S}(M)$
is a Fr\'echet space. We naturally think of elements in $\mathcal{S}(M)$
as real valued functions on $M$ -- by doing so we implicitly choose representatives.

\subsubsection{Remark}Clearly if $M$ is affine, then Definitions \ref{definitions-for-affine}(1) and \ref{def_schwartz_fun_on_definable} coincide.

\

\subsection{Definition}\label{def-temp-dist-on-definable}
Let $M$ be a definable manifold. Define the \emph{space of tempered distributions
on $M$}, denoted by $\mathcal{S}^*(M)$, as the space of continuous linear
functionals on $\mathcal{S}(M)$.

\

The following Lemma \ref{temp-dist-lemma_for_definable} generalizes the observations in \ref{definitions-for-affine}(1,2). Its proof is straightforward,
thus omitted.

\subsection{Lemma}\label{temp-dist-lemma_for_definable}
Let $M$ be a definable manifold and $U\subset M$ an open definable
subset. Then, $\text{Ext}_U^M:\mathcal{S}(U)\hookrightarrow \mathcal{S}(M)$
is a closed embedding, and the restriction map $\mathcal{S}^*(M)\to\mathcal{S}^*(U)$
is onto.

\subsection{Lemma}\label{lemma-to-define-temp-on-non-affine} Let $M$ be a
definable manifold, and let $t:M\to\bR$ be some function. Then, the following
conditions are equivalent:

\

\begin{enumerate}
\item $t\cdot s\in\mathcal{S}(M)$, for any $s\in\mathcal{S}(M)$;
\item $t|_{U}\in\mathcal{T}(U)$, for any affine definable patch $U\subset
M$;
\item there exists an affine definable open cover $M=\bigcup\limits_{i=1}^{m}U_i$
such that for any $1\leq i\leq m$, $t|_{U_i}\in\mathcal{T}(U_i)$.
\end{enumerate}

Proof: Clearly $(2)\Rightarrow(3)$.

$(3)\Rightarrow(1)$: let $s\in\mathcal{S}(M)$, i.e. $s=\sum\limits_{i=1}^{m}\text{Ext}_{U_i}^{M}(s_i)$
for some $s_i\in\mathcal{S}(U_i)$. As $t|_{U_i}\in\mathcal{T}(U_i)$, also
$t|_{U_i}\cdot s_i\in\mathcal{S}(U_i)$. Thus $t\cdot s=\sum\limits_{i=1}^{m}\text{Ext}_{U_i}^{M}(t|_{U_i}\cdot
s_i)\in\mathcal{S}(M)$.

\

$(1)\Rightarrow(2)$: let $U_1:=U\subset M$ be some affine definable patch,
and complete $U_1$ to some affine open definable cover $M=\bigcup\limits_{i=1}^{m}U_i$.
Let $s\in\mathcal{S}(U_1)$. We need to show that $t|_{U_1}\cdot s\in\mathcal{S}(U)$.
Consider $\hat{s}:=\text{Ext}_{U_1}^{M}(s)\in\mathcal{S}(M)$, then $t\cdot
\hat{s}\in\mathcal{S}(M)$. For any $z\in\partial U$ we have by the chain
rule that $t\cdot \hat{s}$ is flat at $z$. Choose an open definable embedding
$\varphi:U\to S^n$ (see the proof of Lemma \ref{tempered_can_separate_disjont_closed}
for further detail). We need to show that for any $y\in\partial(\varphi(U))\subset
S^n$, $(t\cdot \hat{s})\circ \varphi^{-1}$ is flat at $y$, where $\varphi^{-1}$
is the inverse of $\varphi:U\to \varphi(U)$ -- to be more precise we mean
that $\text{Ext}_{\varphi(U)}^{S^n}((t\cdot \hat{s})\circ \varphi^{-1})$
is flat at $y$, but we do not write it to shorten the notation. We take some
sequence $y_1,y_2,y_3,\dots\in \varphi(U)$ converging to $y$. If the sequence
$\{\varphi^{-1}(y_i)\}_{i=1}^{\infty}$ has an accumulation point in $M$,
then by the previous claim $t\cdot \hat{s}$ is flat at this accumulation
point, and $(t\cdot \hat{s})\circ \varphi^{-1}$ is flat at $y$. Otherwise,
we write $t\cdot \hat{s}=\sum\limits_{i=1}^{m}\text{Ext}_{U_i}^{M}(s_i)$,
for some $s_i\in\mathcal{S}(U_i)$. Consider $\text{Ext}_{U_i}^{M}(s_i)$.
As $\{\varphi^{-1}(y_i)\}_{i=1}^{\infty}$ does not have an accumulation point
in $M$, and in particular does not have an accumulation point in $U_i$, we
have that $Ext_{U_i}^{M}(s_i)\circ\varphi^{-1}$ is flat at $y$. We conclude
that also $(t\cdot \hat{s})\circ\varphi^{-1}=\sum\limits_{i=1}^{m}(\text{Ext}_{U_i}^{M}(s_i)\circ\varphi^{-1})$
is flat at $y$. \qed

\subsection{Definition}\label{def_temp_fun_on_definable} Let $M$ be a
definable manifold. A real valued function $t:M\to\bR$ is called \emph{a
tempered function on $M$} if it satisfies the equivalent conditions of Lemma
\ref{lemma-to-define-temp-on-non-affine}. Denote the space of all tempered
functions on $M$ by $\mathcal{T}(M)$.

\subsection{Lemma}\label{lemma-tempered-is-an-algebra-definable}Let $M$
be a definable manifold. Then, $\mathcal{T}(U)$ together with point-wise
addition and multiplication is an $\bR$-algebra. Moreover, any tempered function
that is bounded from below by some positive number is invertible in this
algebra.

\

Proof: It follows easily from \ref{definitions-for-affine}(3) and Definition \ref{def_temp_fun_on_definable}. \qed

\

We now generalize Proposition \ref{tempered-partition-of-unity-affine}:

\subsection{Proposition (tempered partition of unity on definable
manifolds)}\label{tempered-partition-of-unity-definable}Let
$M$ be a definable manifold and $M=\bigcup\limits_{i=1}^{m}U_i$ be a finite
definable open cover. Then, there exist non-negative tempered functions
$\varphi_1,\dots,\varphi_m\in\mathcal{T}(M)$ such that:\begin{enumerate}
                                                       \item $\text{supp}(\varphi_i)\subset
U_i$ for any $i=1,\dots,m$;
                                                       \item $\sum\limits_{i=1}^{m}\varphi_i(x)=1$
for any $x\in M$;
                                                     \end{enumerate}

\

Proof: Following the proof of Proposition \ref{tempered-partition-of-unity},
one easily sees it is enough to generalize Lemma \ref{tempered_can_separate_disjont_closed},
and in fact it is enough to prove the following:

\subsubsection{Lemma}\label{tempered_can_separate_disjont_closed_definable}
Let $A_0,A_1$ be disjoint definable closed subsets of a definable manifold
$M$. Then, there exists a non-negative $t\in\mathcal{T}(M)$ such that $t|_{A_0}=0$
and $t(x)\geq 1$ for any $x\in A_1$.

Proof: Let $M=\bigcup\limits_{i=1}^m X_i$ be an affine definable open
cover. The case $m=1$ follows from Lemma \ref{tempered_can_separate_disjont_closed}.
For simplicity we prove here the case $m=2$. The general case follows in
the same way. Thus we assume $M=X_1\cup X_2$, where $X_1,X_2$ are affine
definable manifolds. As in the proof of Proposition \ref{tempered-partition-of-unity}
we apply the shrinking lemma (that follows from Lemma \ref{normality_of_definable_lemma})
twice to obtain an affine definable open cover $M=Z_1\cup Z_2$ where $\bar{Z_i}\subset
X_i$ for $i=1,2$, and another affine definable open cover $M=Y_1\cup Y_2$
where $\bar{Y_i}\subset Z_i$ for $i=1,2$.

\

Note that $A_0\cap X_1,A_1\cap X_1\subset X_1$ are disjoint closed definable
sets, and so by Proposition \ref{tempered-partition-of-unity-affine} (applied to their complements) there
exists a non-negative $\alpha_1\in \mathcal{T}(X_1)$ such that $\alpha_1|_{A_0\cap
X_1}=0$ and $\alpha_1|_{A_1\cap X_1}=1$. Note that $X_1\setminus Y_1,X_1\setminus(X_1\cap
Y_2)\subset X_1$ are closed definable subsets and, as $M=Y_1\cup Y_2$,
they are disjoint. Thus, there exists a
non-negative $\beta_1\in \mathcal{T}(X_1)$ such that $\beta_1|_{X_1\setminus
Y_1}=0$ and $\beta_1|_{X_1\setminus(X_1\cap Y_2)}=1$. As $\bar{Y_1}\subset
X_1$, we have that $\text{Ext}_{X_1}^{M}(\beta_1)$ is a smooth function on
$M$. Moreover, $X_2\setminus\overline{(X_2\cap Y_1)}$ is affine, $\text{Ext}_{X_1}^{M}(\beta_1)|_{X_2\setminus\overline{(X_2\cap
Y_1)}}=0$, and $M=(X_2\setminus\overline{(X_2\cap Y_1)})\cup X_1$. Thus,
by Lemma \ref{lemma-to-define-temp-on-non-affine}, $\text{Ext}_{X_1}^{M}(\beta_1)\in\mathcal{T}(M)$.
Moreover, similar arguments show that $\text{Ext}_{X_1}^{M}(\beta_1\cdot
t)\in\mathcal{T}(M)$ for any $t\in\mathcal{T}(X_1)$. In particular $t_1:=\text{Ext}_{X_1}^{M}(\beta_1\cdot
\alpha_1)\in\mathcal{T}(M)$. Note that $t_1$ is a non-negative function,
and moreover $t_1|_{A_1\cap (X_1\setminus (X_1\cap Y_2))}=1$ and $t_1|_{A_0}=0$.
Similarly we construct a non-negative $t_2\in\mathcal{T}(M)$ such that $t_2|_{A_1\cap
(X_2\setminus (X_2\cap Y_1))}=1$ and $t_2|_{A_0}=0$.

\

Now consider the affine definable manifold $X:=X_1\cap X_2$, and the disjoint
closed definable subsets $A_0\cap X,A_1\cap X\subset X$. By Proposition \ref{tempered-partition-of-unity-affine} (applied
to their complements)
there exists a
non-negative $\alpha\in \mathcal{T}(X)$ such that $\alpha|_{A_0\cap
X}=0$ and $\alpha|_{A_1\cap X}=1$. Note that $\overline{Y_1\cap Y_2},\overline{(X_1\setminus
Z_1)\cup(X_2\setminus Z_2)}\subset X$ are disjoint closed definable subsets,
so there exists a non-negative $\beta\in \mathcal{T}(X)$ such that $\beta|_{\overline{(X_1\setminus
Z_1)\cup(X_2\setminus Z_2)}}=0$ and $\beta|_{\overline{Y_1\cap Y_2}}=1$.
As above, we conclude that $t_3:=\text{Ext}_{X}^{M}(\beta\cdot\alpha)\in\mathcal{T}(M)$
is a non-negative function, such that $t_3|_{A_1\cap{\overline{Y_1\cap Y_2}}}=1$
and $t_3|_{A_0}=0$.

\

Finally, we define $t:=t_1+t_2+t_3$, and as $(X_1\setminus(X_1\cap Y_2))\cup
(X_2\setminus(X_2\cap Y_1))\cup\overline{Y_1\cap Y_2}=M$, we have that $t$
is a non-negative function such that $t|_{A_0}=0$
and $t(x)\geq 1$ for any $x\in A_1$. This proves Lemma \ref{tempered_can_separate_disjont_closed_definable}
and so Proposition \ref{tempered-partition-of-unity-definable}. \qed

\

Lemmas \ref{lemma_temp_move_s_to_s_definable}
and Proposition \ref{tempered-is-a-sheaf-definable} generalize Lemmas \ref{lemma_temp_move_s_to_s_affine}
and Proposition \ref{tempered-is-a-sheaf-affine} (respectively). They follow
easily from them and from Definition \ref{def_temp_fun_on_definable},
thus we omit their proofs.

\subsection{Lemma}\label{lemma_temp_move_s_to_s_definable}Let $M$ be a
definable manifold, $U\subset M$ a definable open subset, $t\in\mathcal{T}(M)$,
and $s\in\mathcal{S}(M)$. If $\text{supp}(t)\subset U$, then $(t\cdot s)|_{U}\in\mathcal{S}(U)$.

\subsection{Proposition (sheaf property of tempered functions on definable manifolds)}\label{tempered-is-a-sheaf-definable}
Let $M$ be a definable manifold and $M=\bigcup\limits_{i=1}^{m}U_i$ be
a finite definable open cover. Then, the following sequence is exact:

\begin{center}
\leavevmode
\xymatrix{
0 \ar[rr] & & \mathcal{T}(M) \ar[rr]^{\text{Res}_1} & & \bigoplus\limits_{i=1}^{m}\mathcal{T}(U_i)
\ar[rr]^(.4){\text{Res}_2} & & \bigoplus\limits_{i=1}^{m-1}\bigoplus\limits_{j=i+1}^{m}\mathcal{T}(U_i\cap
U_j),}
\end{center}
where the $k$-th coordinate of $\text{Res}_1(t)$ is $\text{Res}_{U_k}^{M}(t)$
and for $1\leq k< l\leq m$ the $(k,l)$-coordinate of $\text{Res}_2(\bigoplus\limits_{i=1}^{m}t_i)$
is $\text{Res}_{U_k\cap U_l}^{U_k}(t_k)-\text{Res}_{U_{k}\cap U_l}^{ U_l}(t_l)$.

\

Finally, we are able to generalize Propositions \ref{tempered-distributions-is-a-sheaf-affine}
and \ref{schwartz_functions-is-a-cosheaf-affine}. The proofs, which we omit,
are the same, where one has to use Lemma \ref{temp-dist-lemma_for_definable},
Proposition \ref{tempered-partition-of-unity-definable}, Lemma \ref{lemma_temp_move_s_to_s_definable}
and Proposition \ref{tempered-is-a-sheaf-definable} instead of observations
\ref{definitions-for-affine}(1,2),
Proposition \ref{tempered-partition-of-unity-affine}, Lemma \ref{lemma_temp_move_s_to_s_affine}
and Proposition \ref{tempered-is-a-sheaf-affine}, respectively:

\subsection{Proposition (sheaf property of tempered distributions on definable manifolds)}\label{tempered-distributions-is-a-sheaf-definable}Let $M$ be a definable manifold and $M=\bigcup\limits_{i=1}^{m}U_i$
be a finite definable open cover. Then, the following sequence is exact:

\begin{center}
\leavevmode
\xymatrix{
0 \ar[rr] & & \mathcal{S}^*(M) \ar[rr]^{\text{Res}_1} & & \bigoplus\limits_{i=1}^{m}\mathcal{S}^*(U_i)
\ar[rr]^(.4){\text{Res}_2} & & \bigoplus\limits_{i=1}^{m-1}\bigoplus\limits_{j=i+1}^{m}\mathcal{S}^*(U_i\cap
U_j),}
\end{center}
where the $k$-th coordinate of $\text{Res}_1(\xi)$ is $\xi|_{\mathcal{S}(U_k)}$
and for $1\leq k< l\leq m$ the $(k,l)$-coordinate of $\text{Res}_2(\bigoplus\limits_{i=1}^{m}\xi_i)$
is $\xi_k|_{\mathcal{S}(U_k\cap U_l)}-\xi_l|_{\mathcal{S}(U_k\cap U_l)}$.

\subsection{Proposition (cosheaf property of Schwartz functions on definable manifolds)}\label{schwartz_functions-is-a-cosheaf-definable}Let $M$ be a definable manifold and $M=\bigcup\limits_{i=1}^{m}U_i$
be a finite definable open cover. Then, the following sequence is exact:

\begin{center}
\leavevmode
\xymatrix{
\bigoplus\limits_{i=1}^{m-1}\bigoplus\limits_{j=i+1}^{m}\mathcal{S}(U_i\cap
U_j) \ar[rr]^(.6){\text{Ext}_1} & & \bigoplus\limits_{i=1}^{m}\mathcal{S}(U_i)
\ar[rr]^{\text{Ext}_2} & & \mathcal{S}(M) \ar[rr] & & 0,}
\end{center}
where the $k$-th coordinate of $\text{Ext}_1(\bigoplus\limits_{i=1}^{m-1}\bigoplus\limits_{j=i+1}^{m}s_{i,j})$
is $\sum\limits_{i=1}^{k-1} \text{Ext}_{U_k\cap U_i}^{U_k}(s_{i,k})-\sum\limits_{i=k+1}^{m}
\text{Ext}_{U_k\cap U_i}^{U_k}(s_{k,i})$, and $\text{Ext}_2(\bigoplus\limits_{i=1}^{m}s_i):=\sum\limits_{i=1}^{m}\text{Ext}^{M}_{U_i}(s_i)$.

\

We end this section by discussing compactly supported functions. First, let
us see what happens on compact manifolds:

\subsection{Proposition}\label{prop_on_compact_any_smooth_is_schwartz}If
$M$ is a compact definable manifold, then $\mathcal{S}(M)=C^\infty(M)$.

\

Proof: It is enough to show that any smooth function is Schwartz. Let $s\in
C^\infty(M),M=\bigcup\limits_{i=1}^{m}U_i$ an affine definable open cover,
and $\varphi_1,\dots,\varphi_m\in\mathcal{T}(M)$ as in Proposition \ref{tempered-partition-of-unity-definable}.
Fix $1\leq j \leq m$. We claim that $s|_{U_j}\cdot \varphi_j|_{U_j}\in\mathcal{S}(U_j)$:
by abuse of notation we think of $U_j$ as an open subset of $S^{n}$ (see
the proof of Lemma \ref{tempered_can_separate_disjont_closed} for further
detail). Clearly $s|_{U_j}\cdot \varphi_j|_{U_j}$ is smooth on $U_j$, so
the only thing to verify is that $\text{Ext}_{U_j}^{S^n}(s|_{U_j}\cdot \varphi_j|_{U_j})$
is smooth on $S^n$. This clearly follows from compactness of $M$ and the
fact that $\varphi_j$ is supported in $U_j$. Thus $s=\sum\limits_{i=1}^{m}\text{Ext}_{U_i}^{M}(s|_{U_i}\cdot
\varphi_i|_{U_i})\in\mathcal{S}(M)$. \qed

\

The following Theorem \ref{prop_any_compact_supp_fun_is_schwartz} is a generalization
of Proposition \ref{prop_on_compact_any_smooth_is_schwartz} above. Its proof
is very similar, thus it is left to the reader.

\subsection{Theorem}\label{prop_any_compact_supp_fun_is_schwartz}Let $M$
be a definable manifold. Then, any smooth compactly supported function
is Schwartz, i.e. $C_c^\infty(M)\subset\mathcal{S}(M)$.

\section{Schwartz functions with full support -- Whitney's theorem in action}\label{sec-whitney}The
goal of this section is to show that on any open subset of $\bR^n$ (not necessarily
definable) there exists a strictly positive Schwartz function (Theorem
\ref{thm-schwartz-with-full-support}), and as a corollary on any definable
manifold there exists a strictly positive Schwartz function (Corollary \ref{thm-schwartz-with-full-support_cor}).
The key ingredient is the following: for a smooth function $f:\bR^n\to\bR$,
it is clear that the flat locus of $f$ -- the set of all points in $\bR^n$
where the Taylor series of $f$ is identically zero -- is a closed subset
of $\bR^n$. We show that the converse is also true: any closed subset of
$\bR^n$ may be realized as the flat locus of some smooth function $f:\bR^n\to\bR$.
We moreover show that this function $f$ can be chosen to be a Schwartz function.

\

We start by recalling a classical theorem of Whitney:

\subsection{Theorem (Whitney)}\label{whitney_zero_locus_theorem}Let $Z\subset\bR^n$
be a closed subset. There exists a smooth function $f:\bR^n\to\bR$ such that
$f^{-1}(\{0\})=Z$.

\

Theorem \ref{flat_locus_theorem} below is a stronger version of Theorem \ref{whitney_zero_locus_theorem}.
The proof of Theorem \ref{flat_locus_theorem} given below is almost
identical to the proof of Theorem \ref{whitney_zero_locus_theorem}, that may be found, for instance, in \cite{Mug}.

\subsection{Theorem}\label{flat_locus_theorem}Let $Z\subset\bR^n$ be a closed
subset. There exists a smooth function $f:\bR^n\to\bR$ such that $f(x)>0$
for any $x\not\in Z$, and $f$ is flat at $Z$.

\

Proof: As $U:=\bR^n\setminus Z$ is open, it can be presented as a countable
union of open balls: $U=\bigcup\limits_{i=0}^{\infty}B_i$, where $B_i=\{x\in\bR^n:|x-c_i|<r_i\}$
is an open ball of radius $r_i>0$ around $c_i$. We choose smooth functions
$\psi_i:\bR^n\to\bR$ such that the following hold:\begin{enumerate}
                                                                        
                     \item $\psi_i(x)>0$ if $x\in B_i$.
                                                                        
                     \item $\psi_i$ is flat at $\bR^n\setminus B_i$, i.e.
for any $x\in \bR^n\setminus B_i$ and any multi-index $\alpha\in(\mathbb{N}_0)^n$,
$\psi_i^{(\alpha)}(x)=0$.
                                                                        
                     \item $\psi_i^{(\alpha)}$ is uniformly bounded by $2^{-i}$
for any multi-index $\alpha\in(\mathbb{N}_0)^n$ satisfying $|\alpha|\leq
i$.
                                                                        
                   \end{enumerate}
                                                                        
                   To show such functions exist we first define smooth functions
$\widetilde\psi_i:\bR^n\to\bR$ by $\widetilde\psi_i(x):=F(\frac{|x-c_i|}{r_i})$,
where $F:\bR\to\bR$ is defined by $F(x):=e^{\frac{1}{x^2-1}}$ if $x\in (-1,1)$,
and $F(x):=0$ otherwise. Clearly $\widetilde\psi_i(x)>0$ if $x\in B_i$ and
$\widetilde\psi_i(x)=0$ otherwise. Moreover $\widetilde\psi_i$ is flat at
$\bR^n\setminus B_i$. Since $\widetilde\psi_i$ and all of its derivatives
$\widetilde\psi_i^{(\alpha)}$ are bounded, we may multiply each $\widetilde\psi_i$
by some sufficiently small positive number $\epsilon_i>0$ to get the desired
$\psi_i's$.

Define $f=\sum\limits_{i=0}^{\infty}\psi_i$. Property (3) of the $\psi_i$'s
imply that $f$ is indeed a smooth function on $\mathbb{R}^n$: $$\abs{\sum\limits_{i=0}^{\infty}\psi_i^{(\alpha)}(x)}\leq
\abs{\sum\limits_{i=0}^{|\alpha|-1}\psi_i^{(\alpha)}(x)}+\abs{\sum\limits_{i=|\alpha|}^{\infty}\psi_i^{(\alpha)}(x)}
\leq \abs{\sum\limits_{i=0}^{|\alpha|-1}\psi_i^{(\alpha)}(x)}+\sum\limits_{i=|\alpha|}^{\infty}2^{-i},$$
and so for any multi-index $\alpha\in(\mathbb{N}_0)^n$ the sum $\sum\limits_{i=0}^{\infty}\psi_i^{(\alpha)}$
converges uniformly on all of $\bR^n$. Thus, $f$ is a smooth function, $f(x)>0$
for any $x\in U$, and $f^{-1}(\{0\})=Z$. This proves Theorem \ref{whitney_zero_locus_theorem}.

\

The fact that $f$ is flat at $Z$ follows from property (2) of the $\psi_i$'s:
by construction $x\in Z$ if and only if for any $i\in\mathbb{N}_0$ we have
$x\in \bR^n\setminus B_i$. Thus, for any $x\in Z$ and any multi-index $\alpha\in(\mathbb{N}_0)^n$
we have $f^{(\alpha)}(x)=\sum\limits_{i=0}^{\infty}\psi_i^{(\alpha)}(x)=0$,
i.e. $f$ is flat at $Z$. \qed

\

\subsubsection{Remark} We constructed a smooth function $f:\bR^n\to \bR$
such that $f(x)>0$ for any $x\not\in Z$ and $f$ is flat at $Z$, i.e. for
any multi-index $\alpha\in(\mathbb{N}_0)^n$ we have $(f^{(\alpha)})^{-1}(\{0\})\supset
Z$. It is true that $(f^{(0)})^{-1}(\{0\})=f^{-1}(\{0\})=Z$, however for
multi-indices $\alpha\in(\mathbb{N}_0)^n$ such that $|\alpha|>0$ it is not
true in general that $(f^{(\alpha)})^{-1}(\{0\})=Z$.

\

We are now ready to prove the key ingredient:

\subsection{Theorem}\label{thm-schwartz-with-full-support}Let $U\subset\bR^n$
be an open subset. There exists $f\in\mathcal{S}(U)$
such that $f(x)>0$ for any $x\in U$.

\

Proof: Let $p:S^n\setminus\{(1,0,\dots,0)\}\to\bR^n$ be the standard rational
stereographic projection (as in \ref{stereo_projection} above). We have that
$\tilde{U}:=p^{-1}(U)\subset S^n$ is an open subset. Thus $S^n\setminus\tilde{U}$
is closed in $S^n$, which is closed in $\bR^{n+1}$, and so $S^n\setminus\tilde{U}$
is closed in $\bR^{n+1}$. By Theorem \ref{flat_locus_theorem} there exists
a smooth function $\tilde{f}:\bR^{n+1}\to \bR$ such that $f(x)>0$ for any
$x\not\in S^n\setminus\tilde{U}$, and $\tilde{f}$ is flat at $S^n\setminus\tilde{U}$.
We note that: \begin{enumerate}
                                                 \item $S^n\setminus\{(1,0,0,\dots,0)\}$
is Zariski open in $S^n$.
                                                 \item $\tilde{f}|_{S^n}$
is flat at $(1,0,0,\dots,0)$ -- here we mean flat in the sense of \cite[Definition
3.15]{ES}, namely it is a restriction of a smooth function on $\bR^{n+1}$
whose Taylor series at $(1,0,0,\dots,0)$ is identically zero.
                                                 \item $p$ is a biregular
isomorphism of affine algebraic varieties.
                                               \end{enumerate}
Combining (1)-(3) above we may use \cite[Theorem 3.23]{ES}, and conclude
that $\tilde{f}\circ p^{-1}\in\mathcal{S}(\bR^n)$. Clearly $f(x)>0$ for any
$x\in U$. Moreover, as the Taylor series of $\tilde{f}$ around any point
in $S^n\setminus\tilde{U}$ is identically zero, the Taylor series of $\tilde{f}\circ
p^{-1}$ around any point in $\bR^n\setminus U$ is identically zero (this
can be alternatively seen by thinking of $S^n$ as a Nash manifold -- then
$p$ is an isomorphism of affine Nash manifolds, and \cite[Theorem 5.4.1]{AG}
may be used instead of \cite[Theorem 3.23]{ES}. In the Nash category the
notion of flatness is understood intuitively as an intrinsic notion, and
the claim that the Taylor series of $\tilde{f}\circ p^{-1}$ around any point
in $\bR^n\setminus U$ is identically zero is clear). We conclude that $f:=(\tilde{f}\circ
p^{-1})|_U$ is a Schwartz function on $U$, that is everywhere positive. \qed

\subsection{Corollary}\label{thm-schwartz-with-full-support_cor}Let $M$ be
a manifold definable in a polynomially bounded o-minimal structure. Then, there exists $s\in\mathcal{S}(M)$ such that
$s(x)>0$ for any $x\in M$.

\

Proof: This is obvious from Definition \ref{def_schwartz_fun_on_definable}
and Theorem \ref{thm-schwartz-with-full-support}. \qed

\section{The exponential case}\label{sec-miller}The theory developed in Sections \ref{sec-embedded}-\ref{sec-general} may
not be generalized to manifolds definable in o-minimal structures that are
not polynomially bounded: in Example \ref{example-of-all-evil} below we show
that if $\mathcal{R}$ is an o-minimal structure that is not polynomially
bounded,
then Corollary \ref{definable-diffeo-induce-iso-of-schwartz} does not hold,
and so the space of Schwartz functions is not well defined: we explicitly
(and quite easily) find two open subsets of the real line $U,V\subset\bR$,
both definable
in $\mathcal{R}$, and a definable diffeomorphism between these
subsets $\nu:U\to
V$, such that $\nu^*|_{\mathcal{S}(V)}\not\subset \mathcal{S}(U)$.

\subsection{Example}\label{example-of-all-evil} Let  $\mathcal{R}$ is an
o-minimal structure
that is not polynomially bounded. By Miller's Dichotomy Theorem (\ref{Miller_thm})
$\mathcal{R}$ is exponential, i.e. the function $x\mapsto e^x$ (from $\bR$
to $\bR$) is definable in $\mathcal{R}$. The map $\text{exp}:\mathbb{R}\xrightarrow{\sim}
\bR_{>0}$ is thus a definable diffeomorphism, whose inverse (which is also
definable in $\mathcal{R}$ by definition) is $\text{log}:\bR_{>0}\xrightarrow{\sim}
\mathbb{R}$. Take a smooth function $f:\mathbb{R}\to\bR$ such that $f(x)=0$
for any $x<0$ and $f(x)=e^{-x}$ for any $x>1$. Such $f$ clearly exists, and
moreover $f\in\mathcal{S}(\mathbb{R})$. Then, for any $y>e$ we have $\text{log}^*(f)(y)=e^{-\text{log}(y)}=\frac{1}{y}$,
and so $\text{log}^*(f)\notin\mathcal{S}(\bR_{>0})$.

\

In Counterexample \ref{schapira_example} below we show that if $\mathcal{R}$ is an o-minimal structure that is not
polynomially
bounded,
then Lemma \ref{tempered_can_separate_disjont_closed} (tempered Urysohn)
does not hold: we explicitly construct an open subset of the plain $U\subset\bR^2$
that is definable
in $\mathcal{R}$, and two disjoint definable closed subsets $A_0,A_1\subset U$, and
show that there
is no $t\in\mathcal{T}(U)$
such that both $t|_{A_0}=0$ and $t|_{A_1}=1$. Thus, even for "embedded" affine manifolds tempered partition
of unity does not exist, and the theory does not generalize. The author thanks Pierre Schapira for this counterexample.

\subsection{Counterexample for tempered Urysohn lemma in the non-polynomially bounded case}\label{schapira_example}Let $\mathcal{R}$ be an o-minimal structure that is not
polynomially
bounded, then by Miller's Dichotomy Theorem (\ref{Miller_thm})
$\mathcal{R}$ is exponential, i.e. the function $x\mapsto e^x$ (from $\bR$
to $\bR$) is definable in $\mathcal{R}$. We
construct two open definable subsets $V'\subset U\subset\bR^2$ and a closed definable subset
$Y\subset V'$, such that there does not exist $t\in\mathcal{T}(U)$ satisfying
both $t|_Y=1$ and $t|_{U\setminus V'}=0$. As $Y,U\setminus V'\subset U$ are
disjoint closed definable subsets, this is a counterexample for a tempered
Urysohn lemma in the non-polynomially bounded case.

Consider $\bR^2$ with coordinates $(x,y)$. Define open definable subsets $V'\subset
V\subset U\subset\bR^2$: $$U=\{(x,y):y>0\},$$ $$V=\{(x,y):y>0\text{ and }\abs{x}<e^{-\frac{1}{y}}\},$$
$$V'=\{(x,y):y>0\text{ and }\abs{x}<\frac{1}{2}\cdot e^{-\frac{1}{y}}\},$$
and also define the closed definable subset $Y\subset U$: $$Y=\{\frac{1}{2}\geq y>0\text{
and }x=0\}.$$ Let $B\in C^\infty(\bR^2)$ be a compactly supported function
such that $B|_{\{x^2+y^2\leq 1\}}=1$ and $B|_{\{x^2+y^2\geq 2\}}=0$. We will
show that tempered Urysohn lemma for non-polynomially
bounded sets does not hold in
general, by showing that if it does then there exists $S\in\mathcal{S}(U)$
such that $S|_Y=1$. This is impossible as then the Taylor series of $\text{Ext}_U^{\bR^2}(S)$
around the origin is not identically zero. Thus we are only interested in
the behaviour of all functions below around the origin, and for that reason
we have the function $B$ appearing in the definitions.

Consider $f\in C^\infty(U)$ defined by $$f(x,y)=e^{\frac{y}{x^2+y^2}}\cdot
B(x,y).$$ Note that $f|_V\in\mathcal{T}(V)$, as for any $(x,y)\in V$ we have
$$\text{dist}((x,y),\bR^2\setminus V)\leq\text{dist}((0,y),\bR^2\setminus
V)\leq e^{-\frac{1}{y}}.$$

Assume there exists $t\in\mathcal{T}(U)$ such that $t|_Y=1$, and $t|_{U\setminus
V'}=0$. In particular $\text{supp}(t)\subset V$.

Consider $s\in \mathcal{S}(U)$ defined by $$s(x,y)=e^{-\frac{1}{y}}\cdot
B(x,y).$$

By Lemma \ref{lemma_temp_move_s_to_s} we have that $(t\cdot
s)|_V\in\mathcal{S}(V)$, and so (by Proposition \ref{prop_temp_is_like_module_proprty_on_schwartz})
$(f\cdot t\cdot s)|_V\in\mathcal{S}(V)$. However $(f\cdot t\cdot s)|_Y=1$,
which is a contradiction.

\begin{appendix}

\section{Applications in representation theory}\label{appendix-rep-theory}

\subsection{Du Cloux's Schwartz induction}Let us recall some classical constructions
in representation theory, and show how to apply our results in order to study
the induced Schwartz representation introduced by du Cloux:

\subsubsection{Induced representations}Consider a group $G$, a subgroup $H\leq
G$ and a finite dimensional representation $(\pi,V)$ of $H$. The induced
representation of $\pi$ from $H$ to $G$ is the vector space $$\{f:G\to V|f(gh)=\pi(h^{-1})f(g)\text{
}\forall g\in G\text{ }\forall h\in H\}\cap C(G,V),$$ where $C(G,V)$ is some
class of functions from $G$ to $V$, depending on the context, and the action
of $G$ is given by $g.f(g')=f(g^{-1}g')$. For instance, for finite groups
$C(G,V)$ is the space of all functions, and if $G$ is locally compact and
$H$ is closed one may either take $L^2(G/H)$ or the space
of all compactly supported functions mod $H$ (the latter is called compact
induction).

\subsubsection{Schwartz induction}Du Cloux developed the theory of Schwarz
functions
on affine Nash manifolds in \cite{dC}, and in particular constructed \emph{Schwartz
induction}. The interested reader is referred to \cite[Section 2.1]{dC} for
full detail, however the general idea is the following: consider the case
where $G$ above is a Nash group (a semi-algebraic Lie group), and $H$ is
a closed subgroup. Define the space $\mathcal{S}(G,V)$ to be the space of
functions from $G$ to $V$ that are coordinate-wise Schwartz functions. Fix
some Haar measure $dh$ on $H$, and consider the map from $\mathcal{S}(G,V)$
to $C^\infty(G,V)$ defined by $\bar{\varphi}(g):=\int\limits_{h\in H}\pi(h)\varphi(gh)dh$
for any $\varphi\in\mathcal{S}(G,V)$. Note that in general $\bar{\varphi}\not\in\mathcal{S}(G,V)$,
e.g. if $G=H$ then $\bar{\varphi}$ is a constant function, thus either zero
or not Schwartz. Now take the class of functions $C(G,V)$ above to be the
subset of $C^\infty(G,V)$ that is the image of the map above, i.e. $C(G,V)=\{\bar{\varphi}|\varphi\in\mathcal{S}(G,V)\}$.

\subsubsection{Induced representations as spaces of sections}When $G$ is
a topological
group and $H$ a closed subgroup, in many cases the induced representation
may be realized as the space of sections of some vector bundle: consider
the action of $G$ on the space $G\times V$ by $g.(g',x)=(gg',x)$, and consider
the equivalence relation on $G\times V$ defined by $(g,x)\sim(gh,\pi(h^{-1})(x))$
for any $h\in H$. This equivalence relation is invariant under the above
action of $G$, and so $G$ acts on $(G\times V)/_\sim$. Projection to the
first
coordinate $p:(G\times V)/_\sim\to G/H$ makes $(G\times V)/_\sim$ a vector
bundle over $G/H$, with typical fiber $V$. For $g\in G$ and $x\in V$, denote
by $[g,x]$ the equivalence class of $(g,x)$ in $(G\times V)/_\sim$. The action
of $G$ on the space of sections of the vector bundle $p:(G\times V)/_\sim\to
G/H$ is the following: consider a section $\phi:G/H\to(G\times V)/_\sim$,
and denote $[g,\phi_g]:=\phi(gH)$. Then, the action of $G$ is given by $(g.\phi)(g'H)=[gg',\phi_{g'}]=[g',\phi_{g^{-1}g'}]$.
It should be stressed that in the above construction we did not specify what
class of sections is considered. This choice depends on the class of functions
$C(G,V)$ above. For instance, one may take smooth sections etc., depending
on the context.

\subsubsection{Understanding Schwartz induction as a space of sections}In
\cite{AGS}
a special case of the above was analyzed, where in particular a one dimensional
representation (a character) was considered. There, $G\times V$ is a
Nash manifold, however the action of $H$ is not a Nash map, and so the quotient
is not a Nash manifold. The authors constructed
an ad-hoc theory of Schwartz spaces suitable for this situation, and thus
were able to realize du Cloux's Schwartz induction as the space of \emph{Schwartz
sections}.

\subsubsection{Schwartz sections of a general definable vector bundle}Our
construction
of Schwartz functions on manifolds that are definable in polynomially bounded o-minimal structures enables a systematic construction of Schwartz sections,
for more general cases. The basic idea is the following: consider a vector
bundle $p:M\to N$, with a finite dimensional typical fiber $V$, where $M$
and $N$ are definable manifolds and $p$ is a definable map, with a
(finite) definable trivialization.
By local triviality one may define Schwarz sections of this bundle: these
are smooth sections that are coordinate-wise Schwartz functions.

\subsubsection{\bf{Example}}\label{example-good-structure}Let $$\bR^\bR_{\text{an}}:=(\bR,<,0,1,+,\cdot,(f),(x^r)_{r\in\bR}),$$ where $f$ ranges over all restricted analytic functions (i.e. ranges over all functions from $\bR^n$ to $\bR$, $n\in\mathbb{N}$, that vanish identically outside $[-1,1]^n$ and whose restrictions to $[-1,1]^n$ are analytic) and all functions of the form $x^r:\bR\to\bR$ defined by $x^r(a):=\begin{cases}
a^r, \ \text{if} \ a>0 \\
0, \ \text{if} \ a\leq0
\end{cases}$. Then, $\bR^\bR_{\text{an}}$ is a polynomially bounded o-minimal structure (see \cite{vdDM,Miller2}), and it is easy to verify that any function on $\bR^\times$ of the form $x\mapsto\abs{x}^r$ (with some fixed $r\in\bR$) is definable in this structure.

Consider a multiplicative character of the group $\bR^\times$ defined by $g\mapsto\abs{g}^r$, where $r\in\bR\setminus\bQ$. This character is not a Nash function, however it is definable in $\bR^\bR_{\text{an}}$. Thus, the corresponding bundle and actions described above are well defined in
our setting (but not in the Nash, i.e. \cite{AG} setting) and Schwartz sections and Schwartz induction may now be defined. 
\subsection{Existence of tempered distributions}A measure on a locally compact
topological space is a linear functional on the space of compactly supported
continuous real-valued functions on the space. Let $G$ be a group acting
on a topological space $X$. In many cases one knows that there exists a "well
behaved" measure $\mu$ on some $G$-orbit $\mathcal{O}\subset X$, and wants
to deduce the existence of a "well behaved" tempered distribution $\eta$
on $X$, that "extends" $\mu$ in some sense. Making this statement precise
one has to make sense of the notion \emph{tempered distribution}, i.e. realize
what is the space of Schwartz functions on $X$. The following case was proven
by Gourevitch, Sahi and Sayag (for the definitions of the terms appearing
in the theorem see \cite{GSS}):

\subsubsection{Theorem (cf. \cite[Theorem 4.3]{GSS})}\label{GSS-theorem}Let
$G$ be a quasi-elementary $\mathbb{R}$-group, and $X$ a smooth quasi-affine
algebraic variety. Assume that there exists a $G$-orbit $\mathcal{O}\subset
X$ that admits a tempered holonomic $G$-equivariant $V^*$-valued measure
$\mu$. Then there exists a generalized $G$-invariant holonomic distribution
$\eta\in\mathcal{S}^*(X,V)$ and a Zariski open $G$-invariant neighborhood
$U$ of $\mathcal{O}$, such that the restriction of $\eta$ to $U$ equals the
extension of $\mu$ to $U$ by zero.

\

In the above Theorem \ref{GSS-theorem} all objects in question (groups, manifolds,
bundles, maps, etc.) are naturally Nash. In this case the theory of Schwartz
spaces is well established, and in particular the notion \emph{tempered distribution}
is well defined. For quasi-projective varieties non-Nash objects appear,
however it seems they are still naturally \emph{definable in some polynomially bounded o-minimal structure}. From this
it follows that the main obstacle to generalize Theorem \ref{GSS-theorem}
to quasi-projective varieties can be removed using our theory (see \cite[Remark
4.6]{GSS}).

\section{A conjecture on $C^\infty$-diffeomorphisms vs. definable $C^\infty$-diffeomorphisms}\label{appendix-conjectures}

In Example \ref{example-of-all-evil}
we saw that $\text{exp}:\bR\xrightarrow{\sim} \bR_{>0}$ is a diffeomorphism
of smooth manifolds, but it does not induce an isomorphism of the corresponding
Schwartz spaces. However, $\bR$ and $\bR_{>0}$ are indeed isomorphic Nash manifolds, i.e. isomorphic in the category of manifolds definable in the polynomially
bounded o-minimal structure $(\bR,<,0,1,+,\cdot)$: the semi-algebraic map $\varphi:\bR\to\bR_{>0}$
defined by $\varphi(x)=x+\sqrt{x^2+1}$ is a diffeomorphism, whose inverse
is given by $\varphi^{-1}(x):=\frac{1-x^2}{2x}$. In particular, the Fr\'echet spaces $\mathcal{S}(\bR)$ and $\mathcal{S}(\bR_{>0})$ are isomorphic. This suggests the following conjecture:

\subsection{Conjecture}\label{conj_on_nash}Let $U,V\subset\bR^n$ be two open
Nash submanifolds. If there exists a $C^\infty$-diffeomorphism $\varphi:U\to V$, then
there exists a Nash diffeomorphism $\phi:U\to V$.

\subsubsection{Remark}For $n=1$ and $n=2$ Conjecture \ref{conj_on_nash} holds.
This follows immediately from \cite[Corollary 3]{Shi1}. The one dimensional
case may be easily seen, as any open semi-algebraic subset of $\mathbb{R}$
is a finite union of intervals (possibly of infinite length).

\

A more general conjecture is the following:

\subsection{Conjecture}\label{conj_on_definable}Let $U,V\subset\bR^n$ be two open
subsets definable in some  o-minimal structure $\mathcal{R}$. If there exists
a $C^\infty$-diffeomorphism $\varphi:U\to V$, then there exists a
$C^\infty$-diffeomorphism $\phi:U\to V$ that is definable in $\mathcal{R}$.

\subsubsection{Remarks}
\begin{enumerate}
\item For $n=1$ Conjecture \ref{conj_on_definable} holds, as any open definable
subset of $\mathbb{R}$ in an o-minimal structure is a finite union of intervals
(possibly of infinite length).
\item Possibly, Conjecture \ref{conj_on_definable} may be proved using cell decomposition (see
\cite[4.2]{vdDM}). Some other results that may be relevant (and are of similar flavour) may be found in \cite{Fi}: for instance, it is shown that a definable $C^\infty$-diffeomorphism exists if and only if a definable $C^1$-diffeomorphism exists. The latter paper contains results that may be used to give an alternative proof of our Lemma \ref{tempered_can_separate_disjont_closed} (tempered Urysohn) -- we chose to give a direct proof using H\"ormander's classical result (Proposition \ref{Hormander's-lemma}).  
\end{enumerate}

\section{Invariants -- a possible application in fractal geometry}\label{appendix-fractals}

The purpose of this appendix is to suggest certain invariants that may be related to the Schwartz theory, and that may be used, for instance, to study fractals. We do this by defining an equivalence relation on the set of all open subsets of $\bR^n$ (for a fixed $n$) and suggesting what equivalent subsets should have in common. We start with the basic definition:

\subsection{Definition}Two open subsets $U,V\subset\bR^n$ are called \emph{Schwartz equivalent} if there exists a $C^\infty$-diffeomorphism $\nu:U\to V$ such that $\nu^*|_{\mathcal{S}(V)}$ is an isomorphism of the
Fr\'echet spaces $\mathcal{S}(V)$ and $\mathcal{S}(U)$.

\subsubsection{Example} By Corollary \ref{definable-diffeo-induce-iso-of-schwartz} if $\nu:U\to V$ is a $C^\infty$-diffeomorphism that is definable in some polynomially bounded o-minimal structure, then $U$ and $V$ are Schwartz equivalent.

\subsection{Example}\label{exmaple-koch-disc}Let $U\subsetneq\bR^2$ be any
non-empty open simply connected subset. By Riemann's open mapping theorem
$U$ is diffeomorphic to the unit disc $\mathbb{D}:=\{(x,y)\in\bR^2:x^2+y^2<1\}$.
We conclude that all open simply connected subsets of the plane (that are
neither empty nor everything) are isomorphic in the category of smooth manifolds
($\bR^2$ is also isomorphic to $\mathbb{D}$ in this category, but this does
not follow from Riemann's open mapping theorem). In particular, the interior
of Koch's snowflake (see construction in \cite[Figure 0.2]{Fa}) is isomorphic
to the unit disc. Its boundary (Koch's snowflake) has Hausdorff dimension
$\text{log}_3(4)$, that is strictly greater than 1 and strictly smaller than
2 (see calculation in \cite[Example 9.5]{Fa}). The boundary of the unit disc
is a smooth manifold of (Hausdorff) dimension 1. Thus, two open diffeomorphic subsets may have boundaries of different Hausdorff dimensions,
and moreover the fact that one boundary is smooth or not says nothing about
the smoothness of the other.

\subsection{Example}\label{exmaple-boundaries-point-circle} Consider $U:=\bR^2\setminus\{(0,0)\}\subset\bR^2$
and $V:=\{(x,y)\in\bR^2:x^2+y^2>1\}\subset\bR^2$. Define a function
$\varphi:U\to V$ by $\varphi(x,y)=(x+\frac{x}{\sqrt{x^2+y^2}},y+\frac{y}{\sqrt{x^2+y^2}})$.
This $\varphi$ is a bijection, whose inverse is given by $\varphi^{-1}(x,y)=(x-\frac{x}{\sqrt{x^2+y^2}},y-\frac{y}{\sqrt{x^2+y^2}})$.
In fact, $U$ and $V$ are both Nash submanifolds of $\bR^2$, and $\varphi$
is a Nash diffeomorphism, and in particular definable in the polynomially bounded o-minimal structure $(\bR,<,0,1,+,\cdot)$. We conclude that $U$ and $V$
are Schwartz equivalent. However, $\partial
U=\{(0,0)\}$ is a single point, whereas $\partial V=\{(x,y)\in\bR^2:x^2+y^2=1\}$
is the unit circle. In particular $\partial U$ is a smooth manifold of dimension
$0$, whereas $\partial V$ is a smooth manifold of dimension $1$. Thus two
Schwartz equivalent subsets may have boundaries of different Hausdorff
dimensions, even if both are smooth.

\

Example \ref{exmaple-boundaries-point-circle} deals with the case where an
open subset of $\bR^n$ has boundary of Hausdorff dimension strictly less
than $n-1$. It is well known that for bounded sets this cannot happen (for
the reader convenience we give here a proof, as suggested in \cite{I}):

\subsection{Proposition}Let $U\subset\bR^n$ be an open bounded subset. Then,
the Hausdorff dimension of $\partial U$ has Hausdorff dimension at least
$n-1$.

\

Proof: Fix an inner product on $\bR^n$, and let $V\subset\bR^n$ be any linear
subspace of codimension 1. Consider the natural projection $p:\bR^n\to V$.
As $p$ does not increase distances it is a Lipschitz transformation, and
so, by \cite[Corollary 2.4(a)]{Fa}, the Hausdorff dimension of $\partial
U$ is at least the Hausdorff dimension of $p(\partial U)$. As $U$ is bounded,
$p(\partial U)\supset p(U)$, and so the Hausdorff dimension of $p(\partial
U)$ is at least the Hausdorff dimension of $p(U)$. As $U$ is open in $\bR^n$,
$p(U)$ is open in $V$, and so has Hausdorff dimension $n-1$. \qed

\subsection{Definition}An open subset of $\bR^n$ is called \emph{irregular}
if its boundary has Hausdorff dimension strictly greater than $n-1$.

\

We now present two extreme situations, one of them might be true. We formulate
these by the following two questions:

\subsection{Question}\label{thm-same-iso-classes}Let $U,V\subset\bR^n$ be
two open subsets, and assume there exists a $C^\infty$-diffeomorphism $\varphi:U\to
V$. Does this imply that there exists a $C^\infty$-diffeomorphism $\nu:U\to V$ such
that $\nu^*|_{\mathcal{S}(V)}$ is an isomorphism of the
Fr\'echet spaces $\mathcal{S}(V)$ and $\mathcal{S}(U)$? In other words, is the usual equivalence in the $C^\infty$-smooth-category equivalent to Schwartz equivalence? 

\subsection{Question}\label{thm-irreg-is-invar}Let $U,V\subset\bR^n$ be two
Schwartz equivalent open subsets. Is it
true that $U$ is irregular if and only if $V$ is irregular?

\

Example \ref{exmaple-koch-disc} above shows that it is not possible for the
answers of both Question \ref{thm-same-iso-classes} and Question \ref{thm-irreg-is-invar}
to be positive. It is not obvious that either must be positive. Note that
proving Conjecture \ref{conj_on_definable} will say nothing about which (if any)
is positive, as a subset is never irregular if it is definable in an o-minimal structure.

\

 If
the answer to Question \ref{thm-irreg-is-invar} is positive, then the Schwartz equivalence classes may be used to determine the Hausdorff dimension of some sets. In
that case maybe even more is true, namely the following:

\subsection{Question}\label{conjecture_theorem_fractal}Let $U,V\subset\bR^n$
be two open Schwartz equivalent subsets and assume
$U$ is irregular. Does this imply that $V$ is irregular and that moreover
the Hausdorff dimension of the boundary of $U$ equals the Hausdorff dimension
of the boundary of $V$?

\

\subsection{Example}\label{example_cantor}Let $U\subset\bR$ be the complement
in $[0,1]$ to the standard $\frac{1}{3}$-Cantor set. As a smooth manifold,
$U$ is nothing more than a disjoint union of countably many open intervals.
Thus, as a smooth manifold, $U$ is isomorphic to the disjoint union of open
intervals of length $\frac{1}{2}$ centered at the integers (denote this union
by $V$). $U$ is clearly irregular, as it has the $\frac{1}{3}$-Cantor set
as its boundary, whereas $V$ has a countable boundary, and so it is not irregular.
The first stage in order to answer the above questions, is to answer the
following question: are $U$ and $V$ Schwartz equivalent?

\

The author conjectures that the answer is negative, and that in general the
answer to Question \ref{conjecture_theorem_fractal} is positive. If this
is not the case and $U$ and $V$ above are Schwartz equivalent, then the answer to Question \ref{thm-irreg-is-invar} is negative (and in
particular the answer to Question \ref{conjecture_theorem_fractal} is also
negative). In that case perhaps a weaker version is true:

\subsection{Question}Let $U,V\subset\bR^n$ be two open \emph{bounded} Schwartz equivalent subsets. Is it true that $U$ is
irregular if and only if $V$ is irregular?

\section{Relation to the Nash and the algebraic categories}\label{appendix-nash}

Schwartz functions on Nash manifolds and on algebraic varieties were defined
and studied in \cite{dC,AG} and in \cite{ES} respectively. The main goal
of this appendix is to show that our theory
is consistent with these works.

\subsection{Compatibility of our theory with known categories}The simplest example of an o-minimal structure is $(\bR,<,0,1,+,\cdot)$. This structure is polynomially bounded, and the sets, morphisms etc. defined in this structure are usually called \emph{semi algebraic}. The category of smooth manifolds definable in this structure is often called the category of ($C^\infty$-)Nash manifolds, or simply the Nash category.

\subsubsection{Affine Nash manifolds}Schwartz functions on Nash submanifolds
of $\bR^n$
were studied extensively in \cite{dC} and in \cite{AG}. In these works an
intrinsic (much more complicated) definition was given (see \cite[Definition
4.1.1]{AG}). However, in the case of open semi-algebraic
subsets of $\bR^n$ this intrinsic definition coincides with the above Definition
\ref{schwartz space for euclidean open definition} -- this is an easy implementation
of \cite[Theorem 5.4.1]{AG}.

\subsubsection{Affine algebraic varieties}Note that an open algebraic subset of
$\bR^n$ (i.e. the non-vanishing locus of a polynomial) is in particular a
Nash submanifold, and an algebraic (biregular) isomorphism of such objects
is a Nash diffeomorphism. Thus, wherever we consider open subsets of $\bR^n$
that are algebraic, our definitions are also compatible with the theory of
Schwartz functions on (affine) algebraic varieties, developed in \cite{ES}.

\subsubsection{General Nash manifolds and general smooth algebraic varieties}It follows
from \cite[Remark I.5.12]{Shi2} that the definitions of Schwartz spaces
introduced above and introduced in \cite{AG} coincide. In particular any
smooth algebraic variety naturally defines a definable manifold, %(then
%the topology has to be refined as explained in \cite[Remark 3.2.15(a)]{BCR})
and again, the definitions of Schwartz spaces introduced above and introduced
in \cite{ES} coincide.

\

The following corollary is a consequence of the compatibility of the theory
introduced in this paper and the theory introduced
in \cite{AG}. It  follows immediately from the fact that any Nash manifold
can be "glued" from finitely many copies of $\bR^n$ by Nash diffeomorphisms
(\cite[Remark I.5.12]{Shi2}) and from the definitions and properties of Nash
functions and Nash differential operators (in particular from the fact that
any Nash function on $\bR^n$ is bounded by a polynomial). For exact definitions
see \cite{AG}.

\subsection{Corollary} Let $M$ be a Nash manifold, $t\in C^\infty (M)$.
Then, the following are equivalent: 

\begin{enumerate}
\item $t$ is a tempered function in the sense of \cite{AG}, i.e. there exists
an (open) affine Nash cover $M=\bigcup\limits_{i=1}^k U_i$ such that for
any $1\leq i\leq k$ the following holds: for any Nash differential operator
$D$
on $U_i$ there exists a Nash function $f$ on $U_i$, such that $\abs{D(t|_{U_i})(x)}\leq
f(x)$ for any $x\in
U_i$.
\item $t\cdot s\in\mathcal{S}(M)$ for any $s\in\mathcal{S}(M)$.
\end{enumerate}

\end{appendix}


\begin{thebibliography}{MMM}

\bibitem[AG08]{AG} A.~Aizenbud, D.~Gourevitch, {\em Schwartz functions on Nash manifolds},
International Mathematics Research Notices (2008), DOI:10.1093/imrn/rnm155.

\bibitem[AGS15]{AGS} A.~Aizenbud, D.~Gourevitch, S.~Sahi, {\em Twisted homology for the mirabolic niradical},
Israel Journal of Mathematics {\bf 206} (2015), 39--88, DOI:10.1007/s11856-014-1150-3.

%\bibitem[BCR98]{BCR} J.~Bochnak, M.~Coste, M-F.~Roy, {\em Real algebraic %geometry},
%Springer-Verlag Berlin Heidelberg (1998), DOI:10.1007/978-3-662-03718-8.

%\bibitem[BM88]{BM1} E.~Bierstone, P.~D.~Milman, {\em Semianalytic and subanalytic
%sets},
%Publications math\'ematiques de l'IHES, {\bf 67} (1988) 5--42.

\bibitem[dC91]{dC} F.~du Cloux, {\em Sur les repr\'esentations diff\'erentiables des groupes de Lie al\'egbriques},
Annales scientifiques de l'\'E.N.S , {\bf $4^e$ s\'erie, tome 24, $n^o$ 3} (1991) 257--318.

\bibitem[vdDM96]{vdDM} L.~van den Dries, C.~Miller, {\em Geometric categories and o-minimal structures},
Duke Math. J. {\bf 84} (1996), 497--540.

\bibitem[ES18]{ES} B.~Elazar, A.~Shaviv, {\em Schwartz functions on real algebraic varieties}, Canadian Journal of Mathematics (to appear), DOI:10.4153/CJM-2017-042-6.

\bibitem[Fa97]{Fa} K.~Falconer, {\em Fractal Geometry -- Mathematical Foundations and Applications},
John Wiley \& Sons Ltd (1997), ISBN:0-471-96777-7.

\bibitem[Fi08]{Fi} A.~Fischer, {\em Smooth functions in o-minimal structures},
Advances in Mathematics (2008), DOI:10.16/j.aim.2008.01.002.

\bibitem[Fr82]{Fr} F.~G.~Friedlander, {\em Introduction to the theory of distributions},
Cambridge University Press (1982), ISBN:0-521-28591-7.

%\bibitem[Ga68]{Ga} A.~M.~Gabrielov, {\em Projections of semi-analytic sets},
%Functional Anal. Appl., {\bf 2} (1968), 282--291 = Funkcional. Anal. i Prilozen. %{\bf 2}, No. 4 (1968), 18--30.


\bibitem[GSS18]{GSS} D.~Gourevitch, S.~Sahi, E.~Sayag, {\em Analytic Continuation of Equivariant Distributions},
International Mathematics Research Notices (2018), DOI:10.1093/imrn/rnx326.

\bibitem[H83]{H} L.~H\"ormander, {\em The Analysis of Linear Partial Differential Operators I},
Springer-Verlag Berlin Heidelberg NewYork Tokyo (1983), ISBN-13: 978-3-642-96752-8.

\bibitem[I10]{I} S.~Ivanov, {\em Hausdorff dimension of the boundary of an open set in the Euclidean space - lower bound},
https://mathoverflow.net/q/40660.

\bibitem[KS96]{KS1} M.~Kashiwara, P.~Schapira, {\em Moderate and formal cohomology associated with constructible sheaves},
M\'emoires Soc. Math. France, {\bf 64} (1996).

\bibitem[KS01]{KS2} M.~Kashiwara, P.~Schapira, {\em Ind-sheaves},
Ast\'erisque, Soc. Math. France, {\bf 271} (2001).

\bibitem[Ma66]{Ma} B.~Malgrange, {\em Ideals of differentiable functions},
Tata Institute of Fundamental Research, Bombay, Oxford University Press (1966).

\bibitem[Mi94a]{Miller} C.~Miller, {\em Exponentiation is hard to avoid},
Proc. Amer. Math. Soc. {\bf 122} (1994), 257--259.

\bibitem[Mi94b]{Miller2} C.~Miller, {\em Expansions of the real field with power functions},
Ann. Pure Appl. Logic {\bf 68} (1994), 79--94.

\bibitem[Mug05]{Mug} M.~Muger, {\em An introduction to differential topology, de Rham theory and Morse theory},
http://www.math.ru.nl/$\sim$mueger/diff\_notes.pdf, 2005.

\bibitem[Mun00]{Mun} J.~R.~Munkres, {\em Topology ($2^{nd}$ edition)},
Prentice Hall, Upper Saddle River, NJ 07458 (2000), ISBN 0-13-181629-2.

\bibitem[PS99]{PS} Y.~Peterzil, C.~Steinhorn, {\em Definable Compactness and Definable Subgroups of O-minimal Groups},
J. London Math. Soc. (2), 59 (1999) 769--786.

\bibitem[Sc51]{schwartz} L.~Schwartz, {\em Th\'eorie des distributions},
l'Institut de Mathématique de l'Université de Strasbourg, nos. 9
and 10 ; Actualit\'es Scientifiques et Industrielles, nos. 1091 and
1122. Vol. I, 1950, 148 pp. Vol. II, 1951, 169 pp.

\bibitem[Sh83]{Shi1} M.~Shiota, {\em Classification of Nash Manifolds},
Annales de l'Institut Fourier, {\bf tome 33, $n^o$ 3} (1983), 209--232.

\bibitem[Sh87]{Shi2} M.~Shiota, {\em Nash manifolds},
Lecture Notes in Mathematics, {\bf 1269} (1987).

%\bibitem[Sh3]{Shi3} M.~Shiota, {\em Geometry of Subanalytic and Semialgebraic
%Sets},
%Progress in Mathematics, {\bf 150} (1997).

\bibitem[T67]{T} F.~Tr\'eves, {\em Topological vector spaces, distributions and kernels},
Academic Press (1967).

%\bibitem[Wa83]{Wa} F.~W.~Warner, {\em Foundations of Differential Manifolds %and Lie Groups},
%Springer-Verlag Berlin Heidelberg GmbH (1983), DOI:10.1007/978-1-4757-1799-0.

\end{thebibliography}
\end{document}